\documentclass[a4paper,parskip=half-,fontsize=11pt]{scrarticle} 
\usepackage[english]{babel}
\usepackage[a4paper, top=0.98in, bottom=0.79in, left=0.98in, right=0.98in]{geometry}
\renewcommand{\and}{, }
\newcommand{\lastand}{, and }

\author{Michael Hanke\lastand Roswitha M\"arz}
\newcommand{\authorlist}{
Michael Hanke\\
Department of Mathematics\\
Royal Institute of Technology\\
100 44 Stockholm, Sweden\\[3mm]
Roswitha M\"arz\\
Institute of Mathematics\\
Humboldt University of Berlin\\
10099 Berlin, Germany
}

\title{Canonical Subspaces of Linear Time-Varying Differential-Algebraic Equations and Their Usefulness for Formulating Accurate Initial Condition}

\usepackage{xfrac}
\usepackage{graphicx}
\usepackage{mathptmx}      
\usepackage{amsmath}
\usepackage{amssymb}
\usepackage{amsfonts}
\usepackage{amsthm}
\usepackage{makeidx}
\usepackage{dcolumn}
\usepackage{enumerate}
\usepackage{xcolor}

\newcommand{\pPi}{\mathnormal{\Pi}}
\newcommand{\im}{\operatorname{im}}

\newcommand{\rank}{\operatorname{rank}}
\newcommand{\diag}{\operatorname{diag}}

\newcommand{\Real}{\mathbb{R}}
\newcommand{\Natu}{\mathbb{N}}

\newtheorem{proposition}{Proposition}[section]
\newtheorem{lemma}[proposition]{Lemma}
\newtheorem{theorem}[proposition]{Theorem}
\newtheorem{corollary}[proposition]{Corollary}

\theoremstyle{definition}
\newtheorem{definition}[proposition]{Definition}
\newtheorem{remark}[proposition]{Remark}

\newtheorem{example}[proposition]{Example}

\begin{document}

\maketitle

\abstract{\small
\begin{center}\textbf{Abstract}\end{center}
Accurate initial conditions have the task of precisely capturing and fixing the free integration constants of the flow considered. This is trivial for regular ordinary differential equations, but a complex problem for differential-algebraic equations (DAEs) because, for the latter, these free constants are hidden in the flow. We deal with linear time-varying DAEs and obtain an accurate initial condition by means of applying both a reduction technique and a projector based analysis. The highlighting of two canonical subspaces plays a special role. In order to be able to apply different DAE concepts simultaneously, we first show that the very different looking rank conditions on which the regularity notions of the different concepts (elimination of unknowns, reduction, dissection, strangeness, and tractability) are based are de facto consistent. This allows an understanding of regularity independent of the methods.
}

\section{Introduction}
\label{introduction}
In the context of generalized eigenvalue problems, regular matrix pencils, and descriptor systems, so-called spectral projectors and their associated subspaces receive special attention.
In contrast, they have received little attention in the context of the analysis of differential-algebraic equations (DAEs).
We dedicate the present work to this question.  

In doing so, we take as a basis the \emph{regularity notion} that 
adopts the reduction idea presented in \cite{RaRh}. We will prove that  this notion captures equally
proper
the very different approaches from \cite{Cis1982,Jansen2014,RaRh,KM2006,CRR}. In particular, we provide a set of characteristic values from which one can derive the characteristic values of the different approaches, and which, in this sense, are method-independent.

We deal both with DAEs in standard form,
\begin{align}\label{sDAE}
 E(t)x'(t)+F(t)x(t)=q(t),\quad t\in \mathcal I,
\end{align}
and DAEs with properly involved derivative term,
\begin{align}\label{pDAE}
 A(t)(Dx)'(t)+B(t)x(t)=q(t),\quad t\in \mathcal I,
\end{align}
with sufficiently smooth coefficient functions $E,F:\mathcal I\rightarrow \Real^{m\times m}$,
and $A:\mathcal I\rightarrow \Real^{m\times k}$, $D:\mathcal I\rightarrow \Real^{k\times m}$, $B:\mathcal I\rightarrow \Real^{m\times m}$, respectively.

The properly involved derivative is merely useful when it comes to exact smoothness issues of the solutions. Such properties do not matter here, since 
we assume continuously differentiable solutions, as is widely customary and even necessary for the use of the related concepts.
Here we are interested in the time-varying subspaces of $\Real^m$ where the solution values are located.
The canonical subspaces $S_{can}$, $N_{can}$, as well as the associated projector function $\pPi_{can}$ should not depend on how the leading term of the DAE is written down.

We will prove that any regular DAE \eqref{sDAE} or \eqref{pDAE} with sufficiently smooth coefficients possesses a well-defined canonical projector function $\pPi_{can}$ with associated subspaces $S_{can}$ and $N_{can}$ in analogy to, and as a generalization of, the spectral projector. Then the reduction procedure from \cite{RaRh} provides a basis of the canonical subspace $S_{can}$.

The knowledge of the subspace $N_{can}$ allows to elaborate \emph{accurately stated initial conditions} in the sense of \cite[Definition 2.3]{LMW} and \cite[Theorem 2.52]{CRR}, i.e., the accurate capture of the implicit free integration constants. Using a basis of the corresponding canonical subspace $S_{*\,can}$ to the adjoint DAE, we obtain a suitable matrix which allows to formulate accurate initial conditions.

Accurately stated initial conditions are something quite different from consistent initial values. The latter are needed in integration codes so far, but the generation of  consistent initial data is a problem on its own, e.g., \cite{EsLa2018,BurgerGerdts,BCP89}. 

The least-squares collocation procedures for initial and boundary value problems with higher index differential-algebraic equations proposed and analyzed in \cite{HMTWW,HMT,HM,HM1,HM2021A,HM2021B,Hanke2022} work incredibly good, but one has to provide accurate initial and boundary conditions.
Moreover, for the time-stepping version \cite{HM2020} for systems with dynamic degrees of freedom, one needs also accurately stated transition conditions from one window to the next. 

The paper is structured as follows:
In Section \ref{s.review}, as an introduction, we recapitulate the state of affairs for matrix pencils and DAEs with constant coefficients as well as for semi-explicit index-one DAEs. Then, in Section \ref{s.regular}, we use the reduction approach of \cite{RaRh} , agree on the regularity notion, and provide the first canonical subspace $S_{can}$, which is related to the flow of the homogeneous DAE. The following Section \ref{s.comparison} is concerned with a comparison of different concepts, which in the end will show the general validity of the agreed notion of regularity together with its characteristics. This allows the simultaneous use of the results of the different associated theories. In Section \ref{s.accurateIC}  we elaborate what is meant by accurate initial conditions and the role of the second canonical subspace $N_{can}$. The practically most important result is then stated in Theorem \ref{t.ker} in Section \ref{subs.ker}. In this theorem, we derive a matrix, suitable for formulating accurate initial conditions, by using information coming from the adjoint DAE. Further, a nontrivial example is discussed to some extend in Section \ref{s.CamobellMoore}. There is an appendix with additional information about some technical details.

\section{A review of two well-known cases for motivation, and orientation}\label{s.review}
\subsection{A look back at the well-understood regular matrix pencils}\label{s.pencil}
Given the ordered pair $\{E,F\}$ of matrices $E,F\in\Real^{m\times m}$, we consider the expression  $p(\lambda)=\det(\lambda E+F)$ being a polynomial in $\lambda$, and the so-called matrix pencil $\lambda E+F$. If the polynomial does not vanish identically, the matrix pencil is said to be \emph{regular}. In turn, the pair $\{E,F\}$ is said to be \emph{regular} then.

For each regular pair  $\{E,F\}$, there exist two integers $d,\mu\geq 0$, $\mu\leq m-d$, and nonsingular matrices $L,K\in\Real^{m\times m}$ such that
\begin{align}\label{MP1}
 LEK=\begin{bmatrix}
   I_{d}&\\
   &\mathcal N
  \end{bmatrix},\quad 
  LFK=\begin{bmatrix}
   \mathcal 
   W&\\
   &I_{m-d}
  \end{bmatrix}.
\end{align}
Here and in the following, $I_k\in\Real^{k\times k}$ denotes the identity matrix for $k\in\Natu$.
Thereby, $\mathcal N$ is absent if $d=m$. Otherwise, $\mathcal N$ ist nilpotent with index $\mu$, i.e., $\mathcal N^{\mu}=0$, $\mathcal N^{\mu-1}\neq 0$.
Both integers $\mu, d$, the sizes of the blocks $\mathcal N$ and $\mathcal W$, their eigenstructure as well as the subspaces 
\begin{align}\label{MP2}
 S_{can} :=\im K\begin{bmatrix}
       I_{d}\\0
      \end{bmatrix}\subset\Real^{m},\quad 
N_{can} :=\im K\begin{bmatrix}
       0&\\I_{m-d}
      \end{bmatrix}\subset\Real^{m},
\end{align}
are uniquely determined by the matrix pencil or the pair $\{E,F\}$. Formula \eqref{MP1} describes the \emph{Weierstra{\ss}--Kronecker  form}\footnote{In \cite{BeIlTr12,Trenn}, the term \emph{Quasi-Weierstra{\ss}} is used to emphasize that, unlike the standard Kronecker form, no Jordan form is required for $\mathcal N$ and $\mathcal W$ here.} of the matrix pencil and pair, respectively, and $\mu$ is the \emph{Kronecker index}.
We refer to \cite[Proposition 1.3]{CRR} for corresponding proofs, because there the proof is arranged in such a way that also the less often recognized fact \eqref{MP2} results immediately\footnote{In terms of the proof in \cite[Page 6]{CRR}: Using $\tilde K^{-1}K=\diag (V_{11}, V_{22})$ it simply results that
\begin{align*}
 \im K\begin{bmatrix}
       I_{d}\\0
      \end{bmatrix}
&=\im \tilde K\tilde K^{-1}K\begin{bmatrix}
       I_{d}\\0
      \end{bmatrix}
      =\im \tilde K\begin{bmatrix}
       I_{d}\\0
      \end{bmatrix}V_{11}
      =\im \tilde K\begin{bmatrix}
       I_{d}\\0
      \end{bmatrix},\\
 \im K\begin{bmatrix}
       0\\I_{m-d}
      \end{bmatrix}
&=\im \tilde K\tilde K^{-1}K\begin{bmatrix}
      0\\ I_{m-d}
      \end{bmatrix}
      =\im \tilde K\begin{bmatrix}
      0\\ I_{m-d}
      \end{bmatrix}V_{22}
      =\im \tilde K\begin{bmatrix}
      0\\ I_{m-d}
      \end{bmatrix}.
\end{align*}
}. The \emph{canonical subspaces} $S_{can}$ and $N_{can}$ associated to  the pencil $\lambda E+F$ or pair $\{E,F\}$ have dimensions $d$ and $m-d$, respectively. They intersect trivially only such that the decomposition
\begin{align*}
 S_{can}\oplus N_{can}=\Real^m
\end{align*}
becomes true. The so-called \emph{canonical} or  \emph{spectral projector}\footnote{Also \emph{consistency projector} in \cite{Trenn}.} $\pPi_{can}$ onto $S_{can}$ along $N_{can}$ turns out to be also uniquely determined by $\{E,F\}$ and a possible representation reads
\begin{align*}
 \pPi_{can}=K\begin{bmatrix}
       I_{d}&0\\0&0
      \end{bmatrix}K^{-1}.
\end{align*}
At this place it should be noted that, in the linear algebra context, for a regular pair $\{E,F\}$ with $d< m$, so-called complementary deflating subspaces corresponding to its finite and infinite eigenvalues play their role. $\pPi_{can}$ is then the (right) spectral projector onto the right deflating subspace and $L^{-1}\begin{bmatrix}
       I_{d}&0\\0&0
      \end{bmatrix}L$ is the (left) spectral projector onto the left deflating subspace. The latter hardly plays a role in DAE theory, which is why the addition "right" is omitted here for the spectral projector $\pPi_{can}$.

It is worth adding that always $\ker E\subseteq N_{can}$ holds true.
\medskip

The linear constant coefficient DAE associated with the regular pair $\{E,F\}$,
\begin{align}\label{MP3}
 Ex'(t)+Fx(t)=q(t),\quad t\in \mathcal I,
\end{align}
decomposes corresondingly into
\begin{align*}
 u'(t)+\mathcal Wu(t)&=f(t),\\
 \mathcal Nv'(t)+v(t)&=g(t),\quad t\in\mathcal I,
 \end{align*}
 with
 \begin{align*}
       Lq(t)=\begin{bmatrix}
        f(t)\\g(t)
       \end{bmatrix},\quad
        x(t)=K\begin{bmatrix}
        u(t)\\v(t)
       \end{bmatrix}=\underbrace{K\begin{bmatrix}
        u(t)\\0
       \end{bmatrix}}_{\in S_{can}}+\underbrace{K\begin{bmatrix}
        0\\v(t)
       \end{bmatrix}}_{\in N_{can}}.
\end{align*}
Each solution of the DAE \eqref{MP3} has the form
\begin{align}\label{MP7}
x(t)&=K\begin{bmatrix}
        u(t)\\v(t)
       \end{bmatrix},\\
 u(t)&=\exp^{-(t-a)\mathcal W}u(a)+\int_{a}^{t}\exp^{(s-t)\mathcal W}f(s)\textrm{ds},\nonumber\\ 
 v(t)&=\sum_{j=0}^{\mu-1}(-1)^j \mathcal N^jg^{(j)}(t),\quad t\in\mathcal I,\quad a\in\mathcal I.\label{6a}
 \end{align}
Clearly, the integer $d$ is at the same time the dynamical degree of freedom of the  system \eqref{MP3}. The case $d=0$, and hence $\pPi_{can}=0$, $S_{can}=\{0\}$  may happen. 

At this place, we  emphasize that, when dealing with DAEs, one has to assume that neither the transformation matrices $K,L$ nor the canonical subspaces together with spectral projectors are practically available with little effort. Calculating them  is considerably easier\footnote{Here \emph{Wong sequences} \cite{BeIlTr12,Trenn,Wong74} offer themselves. Furthermore, \cite[Proposition 1.3]{CRR}  is  constructive and the canonical projector can be determined along the line of \cite{Mae1996,Wong2009} and \cite[Part III]{CRR}.} than solving general eigenvalue problems, but still  quite complex and rather hopeless for the case of variable coefficients, which we are actually interested in.

In case of the homogeneous DAE \eqref{MP3} with $d>0$, the flow  is completely  restricted to remain in the subspace $S_{can}$ and in turn, at each time $a\in \mathcal I$, and arbitrary $x_{a}\in S_{can}$, there is exactly one solution passing through, i.e., $x(a)=x_{a}$.
The situation is quite different for nonhomogeneous DAEs. Only if the function $g$ vanishes identically,\footnote{This means, that all $q(t)$ are completely restricted to remain within the left deflating subspace.} and in turn $v(t)\equiv 0$, then the flow is again restricted to $S_{can}$. The property $g(t)\equiv 0$ is practically very unlikely and unrecognizable. Otherwise, one is confronted with nontrivial parts $v(t)$ in \eqref{MP7}. 
In particular, for given $a\in \mathcal I$, $x_a\in S_{can}$, it results that 
\begin{align*}
x_{a}=K\begin{bmatrix}
        u_a\\0
       \end{bmatrix},\quad
 x(a)-x_{a}=\underbrace{K\begin{bmatrix}
        u(a)-u_{a}\\0
       \end{bmatrix}}_{\in S_{can}}+\underbrace{K\begin{bmatrix}
        0\\v(a)
       \end{bmatrix}}_{\in N_{can}},
\end{align*}
and, obviously, the conditions 
\begin{align}
 x(a)-x_{a}\in N_{can},\label{MP8}\\
 \text{or, equivalently, }\quad \pPi_{can}x(a)=x_{a},\nonumber
\end{align}
imply $u(a)=u_{a}$, and vice versa.

Also for $x_a\in \Real^m$ (instead of  $x_{a}\in S_{can}$) we have
\begin{align*}
x(a)-x_a=\underbrace{K\begin{bmatrix}
       u(a)- u_a\\0
       \end{bmatrix}}_{\in S_{can}}+\underbrace{K\begin{bmatrix}
        0\\v(a)-v_a
       \end{bmatrix}}_{\in N_{can}},
\end{align*}
and the conditions 
\begin{align}
 x(a)-x_{a}\in N_{can},\label{MP9}\\
 \text{or, equivalently, }\quad \pPi_{can}x(a)=\pPi_{can}x_{a},\label{MP10}
\end{align}
imply $u(a)=u_{a}$, and vice versa.
We emphasize the unique role of the subspace $N_{can}$ in the DAE context shown here. Only the knowledge of this subspace allows an accurate fixation of the corresponding single solution with property \eqref{MP10}, and this even without knowledge of the projector $\pPi_{can}$ itself.
\medskip

Since each function value $x(a)$ may incorporate a nontrivial $v(a)$, the initial condition $x(a)=x_a\in \Real^m$  includes the necessary condition $v(a)=v_a$ which means that the corresponding differentiations in \eqref{6a} for providing $v(a)$ must be carried out. To really calculate such a kind of \emph{consistent initial values} is quite laborious and uncertain.
However, if $x_{a}$ is not consistent, then the initial value problem (IVP) is not solvable.

On the other hand, if one assumes that the initial conditions are to fix a unique solution from the flow, then exactly $d$
initial conditions shall be set. The following three equivalent conditions are useful,
\begin{align}
 \pPi(x(a)-x_a)&=0,\label{MP4}\\
 x(a)-x_a&\in N_{can},\label{MP5}\\
 G_a (x(a)-x_a)&=0,\label{MP6}
\end{align}
in which the matrix $G_a\in\Real^{d\times m}$ has full row-rank $d$ and its nullspace coincides with $N_{can}$, that is, $\ker G_a=N_{can}$, and $\pPi\in\Real^{m\times m}$ is any projector matrix such that $\ker\pPi=N_{can}$, among others $\pPi_{can}$ is allowed.
Also, for any sufficiently smooth inhomogeneity $q$, an IVP with any of these conditions is uniquely solvable without any further  consistency conditions related to $q$ having to be satisfied. In this sense, these initial conditions are accurately formulated.

We conclude our brief review by showing that various fundamental matrices can be defined for the DAE, all of which are now pointwise singular, have constant rank $d$, and whose pointwise image is exactly $S_{can}$. Typical here are on the one hand the maximal-size version $X(t)$ using the canonical projector, and on the other hand the minimal-size version $X_{b}$ using a basis $C$ of $S_{can}$ given by
\begin{align*}
 EX'(t)+FX(t)&=0,\quad X(a)=\pPi_{can},\\
 EX_{b}'(t)+FX_{b}(t)&=0,\quad X_{b}(a)=C.
\end{align*}
Both versions have their justification. The canonical projector contains more information, namely that about both $N_{can}$ and $S_{can}$, and the corresponding maximal fundamental solution has semigroup properties analogous to the ODE case, e.g., \cite{CRR}, which is very helpful for analytical studies. 
On the other hand, it is usually easier to compute a basis numerically than a projector, especially if it is not an orthoprojector and thus one needs, say, the bases of two subspaces.

Among other things, we also show below a new way to compute the canonical projector $\pPi_{can}$ which seems to be more convenient than the iterative original version in \cite{Mae1996}. In the present paper, by means of the reduction technique presented in Section \ref{s.regular}, a basis $C$ of $S_{can}$ and a basis $C_*$ of the subspace $S_{*\,can}$ associated with the adjoint matrix pair $\{-E^*,F^*\}$ are determined. From  the relations
\begin{align*}
 \im C=S_{can},\quad \ker C^*_{*}E=N_{can}
\end{align*}
which are valid according to Theorem \ref{t.ker}, one can then derive  the projector by usual numerical algebra techniques. Moreover, we may choose $G_a=C^*_{*}E$ for IVPs.

\subsection{Semi-explicit index-one DAEs}\label{s.semi}
Let the matrix function $F:\mathcal I\rightarrow\Real^{m\times m}$ be sufficiently smooth.
The pair given by 
 \begin{align*}
 E=\begin{bmatrix}
    I_r& 0\\0&0
   \end{bmatrix},\quad
F= \begin{bmatrix}
    F_{11}& F_{12}\\F_{21}&F_{22}
   \end{bmatrix},
 \end{align*}
 is associated with the DAE $Ex'+Fx=q$; in more detail,
\begin{align*}
 x_1'(t)+F_{11}(t)x_1(t)+F_{12}(t)x_2(t)&=q_1(t),\\
 F_{21}(t)x_1(t)+F_{22}(t)x_2(t)&=q_2(t).
\end{align*} 
Let the entry $F_{22}(t)\in \Real^{(m-r)\times(m-r)}$ remain nonsingular on $\mathcal I$. Then the DAE 
is regular with index one and posesses the basic subspaces  
\begin{align}
      N&=\ker E=\im \begin{bmatrix}
             0\\I_{m-r}
            \end{bmatrix},\nonumber\\
S(t)&=\{z\in \Real^m: F(t)z\in \im E\}=\im C(t),\quad C(t)=\begin{bmatrix}
             I_r\\-F_{22}(t)^{-1}F_{21}(t)
            \end{bmatrix},\nonumber\\
            &N\cap S(t)=\{0\},\quad t\in\mathcal I.\label{semi1}
     \end{align}
such that $\Real^m=S(t)\oplus N$. 
\medskip

For better clarity, we usually omit the $t$-argument in the following. The relations are then meant pointwise. 

Each solution $x:\mathcal I\rightarrow\Real^m$ of the DAE has the form 
\begin{align}\label{semi2}
x&=\begin{bmatrix}
   Ux_1(a)+f\\-F_{22}^{-1}F_{21}(Ux_1(a)+f)+F_{22}^{-1}q_2
  \end{bmatrix}=
  C(Ux_1(a)+f)+\begin{bmatrix}
                0\\ F_{22}^{-1}q_2
               \end{bmatrix},
  \end{align}
in which $U$ denotes the solution of
\begin{align*}
 U'+(F_{11}-F_{12}F_{22}^{-1}F_{21})U=0,\quad U(a)=I_r,
\end{align*}
 $x_1(a)$ can be predefined arbitrarely by initial conditions, and 
\begin{align*}
 f(t)=U(t)\int_{a}^{t}U(s)^{-1}(q_1(s)-F_{12}(s)F_{22}(s)^{-1}q_2(s))\textrm{ds},\quad t\in\mathcal I.
\end{align*}
Hereby $d=r$ is the dynamical degree of freedom. Further, if $q=0$, in turn $f=0$,  the solutions simplify to
\begin{align*}
x&=\begin{bmatrix}
   Ux_1(a)\\-F_{22}^{-1}F_{21}Ux_1(a)
  \end{bmatrix}=CUx_1(a),\quad x(a)=C(a)x_1(a).
  \end{align*}
Since $C(a)$ is a basis of $S(a)$, and if $S_{can}(a)$ is to be again the subspace containing all solution values of the homogeneous DAE at time $a$, it results that $S=S_{can}$.

Next we are looking for an appropriate complementary to $S_{can}$ subspace $N_{can}$ in view of accurately stated initial conditions. 
The sought subspace must have dimension $m-r$.

Let $x_a=\begin{bmatrix}x_{a1}\\x_{a2}\end{bmatrix}\in\Real^m$ with $x_{a1}\in\Real^r$.
Obviously, the condition $x(a)-x_{a}\in N$ implies $x_1(a)=x_{a\, 1}$, and vice versa. 
Owing to the property \eqref{semi1}, the subspace $N$ is complementary to $S(t)$ for all $t$. 
Here are the projector matrices $P_S(t)$ onto $S(t)$ along $N$ and $P$ onto $N^{\perp}$ along $N$:
\begin{align*}
 P_S=\begin{bmatrix}
             I_r&0\\-F_{22}^{-1}F_{21}&0
            \end{bmatrix},\quad
P=\begin{bmatrix}
             I_r&0\\0&0
            \end{bmatrix}.
\end{align*}
All IVPs
\begin{align*}
 Ex'+Fx=q,\quad x(a)-x_a\in N,
\end{align*}
are uniquely solvable, and each solution satisfies  $x_1(a)=x_{a\, 1}$, the subspace $N$ plays the same role as $N_{can}$ in  Subsection \ref{s.pencil}, so that $N_{can}=N$, $\pPi_{can}=P_{S}$.
\medskip

We underline the uniqueness of $N_{can}$, which means that the condition $x(a)-x_a\in \tilde N$ with any $(m-r)$-dimensional subspace $\tilde N$ other than $N$ involves terms $x_2(a)-x_{a\,2}$, thus parts of $q(a)$.
In particular, the condition $x(a)-x_a\in S(a)^{\perp}$, or equivalently $C(a)^*(x(a)-x_a)=0$, yields\footnote{With an orthonormal basis $C$ for $S$ and the ansatz $x=Cz+(I-CC^*)x$ one can obtain a regular ODE for $z$ together with uniquely solvable IVPs \cite[Remarr 5.2]{Ma2017}. The equation declared as \emph{essential underlying implicit ODE} in \cite{LiMe2012}, for example, is based on such an approach. Then, only as far as homogeneous DAEs are concerned, the initial condition $C(a)^*(x(a)-x_a)=0$ makes sense.}
\begin{align*}
 C(a)^*(x(a)-x_a)=x_1(a)-x_{a\, 1}-(F_{22}(a)^{-1}F_{21}(a))^*(x_2(a)-x_{a\, 2})=0.
\end{align*}
This will be an accurate initial condition, only if $F_{21}=0$, but then $S(a)^{\perp}$ coincides with $N$.
\medskip

It is  evident that accurate initial conditions should precisely determine the value $x_1(a)$, what can be done by one of the conditions
\begin{align}
 x(a)-x_a\in N=N_{can}(a),\label{V1}\\
 P(x(a)-x_{a})=0,\label{V2}\\
 G_a(x(a)-x_{a})=0,\label{V3}
\end{align}
with a matrix $G_a$ such that $\ker G_a=N_{can}$, $\rank G_a=d$.
For arbitrary given $x_a\in \Real^m$, each of these condition yields $x_1(a)=x_{a\,1}$. And this underlines the role of the canonical subspace $N=N_{can}$.
\begin{remark}\label{r.Glattheit}
Like \cite{RaRh}, most work on DAEs assumes continuously differentiable solutions. 
We will do the same in the present paper in order to include all relevant approaches, in particular the one of \cite{RaRh} which is especially easy to understand. 
This means that we must assume somewhat higher smoothness of the data than in those concepts like \cite{CRR} where sharp solvability statements with the lowest possible smoothness properties are required. 
The difference is already visible in the solution representation \eqref{semi2}: To have $x_2$ not only continuous but continuously differentiable, $F^{-1}_{22}F_{21}$ and $F^{-1}_{22}q$ must be continuously differentiable. In concepts with properly involved derivative, one also accepts solution components $x_2$ being only contimuous, thus merely continuous data.

We point out that the matrix $C(t)$ serves as a basis of the subspace $S(t)$.  
The assumption of a continuously differentiable solution goes along with the assumption that $S$ is not only a continuous but also a continuously differentiable subspace varying in $R^m$. Such smooth bases form the starting point on each level of the reduction procedure in \cite{RaRh}.

Taking into account the possible lower smoothnesses in the reduction steps, for example as in \cite{Jansen2014}, requires enormous technical effort that is hiding the essential principles.
\end{remark}

\section{Regular time-varying DAEs and their canonical subspace $S_{can}$
}\label{s.regular}
We turn to the ordered pair $\{E,F\}$ of matrix functions $E,F:\mathcal I\rightarrow\Real^{m\times m}$ being sufficiently smooth, at least continuous, and consider the associated DAE
\begin{align}\label{DAE0}
 E(t)x'(t)+F(t)x(t)=q(t),\quad t\in \mathcal I,
\end{align}
as well as the accompanying time-varying subspaces in $\Real^{m}$,
\begin{align}\label{sub}
 N(t)=\ker E(t),\quad S(t)=\{z\in \Real^m:F(t)z\in\im E(t)\},\quad t\in\mathcal I.
\end{align}
In accordance with Section \ref{s.review} we denote the subspace containing the flow of the homogeneous DAE at time $\bar t$ by $S_{can}(\bar t)$, that is, the set of all possible function values $x(\bar t)$ of solutions of the DAE $Ex'+Fx=0$,
\begin{align*}
 S_{can}(\bar t):=\{\bar x\in \Real^m: \text{there is a solution \;} x:(\bar t-\delta,\bar t+\delta)\cap\mathcal I\rightarrow \Real^m\;\\
 \text{ of the homogeneous DAE such that } x(\bar t)=\bar x\},\quad \bar t\in\mathcal I.
\end{align*}
The second canonical subspace $N_{can}$ is, if it exists, by definition, a pointwise  complement to $S_{can}$, 
\begin{align*}
 \Real^m=N_{can}(\bar t)\oplus S_{can}(\bar t), \quad \bar t\in\mathcal I,
\end{align*}
such that  each IVP,
\begin{align}\label{IVP}
 Ex'+Fx=q,\quad x(\bar t)-\bar x\in N_{can}(\bar t),
\end{align}
with $\bar t\in \mathcal I,$ $\bar x\in \Real^m$, and sufficiently smooth $q$, is uniquely solvable without any consistency conditions for $q$ or its derivatives at the point $\bar t$. This  subspace will be further specified in Section \ref{s.accurateIC} below.

In the present section we 
will agree on what \emph{regular} DAEs are, and show that then the time-varying subspace $S_{can}(\bar t)$ is well-defined on all $\mathcal I$, and has constant dimension.

In Section \ref{s.accurateIC} we will see that in case of  a regular DAE
both canonical subspaces are well-defined with dimensions independent of $\bar t\in \mathcal I$. The associated projector function thus becomes a generalization of the spectral projector for regular matrix pencils in Subsection \ref{s.pencil}.

\begin{definition}\label{d.prereg}
 The pair  $\{E,F\}$ and the DAE \eqref{DAE0}, respectively, are called \emph{pre-regular} on $\mathcal I$ if 
 \begin{align*}
  \im [E(t) \,F(t)]=\Real^m,\quad \rank E(t)=r,\quad \dim N(t)\cap S(t)=\theta, \quad t\in\mathcal I,
 \end{align*}
with integers $0\leq r\leq m$ and $\theta\geq 0$.
Additionally, if $\theta =0$ and $r<m$, then the DAE is called \emph{regular with index one}, but if $\theta =0$ and $r=m$, then the DAE is called \emph{regular with index zero}.
\end{definition}
We underline that any pre-regular pair $\{E,F\}$ features three subspaces $S(t)$, $N(t)$,  and $N(t)\cap S(t)$ having constant dimensions $r$, $m-r$, and $\theta$, respectively. 
\medskip

We emphasize and keep in mind that now not only the coefficients are time dependent, but also the resulting subspaces. Nevertheless, we suppress in the following mostly the argument $t$, for the sake of better readable formulas. The equations and relations are then meant pointwise for all arguments.
\medskip

The different cases for $\theta=0$ are well-understood. 
A regular index-zero DAE is actually a regular implicit  ODE and $S_{can}=S=\Real^m,\, N=\{0\}$. Regular index-one DAEs feature $S_{can}=S,\, N_{can}=N$, e.g., \cite{GM86,CRR}, also Section \ref{s.semi}. Note that $r=0$ leads to $S_{can}=\{0\}$.
All these cases are only interesting here as intermediate results.
\bigskip

We turn back to the general case and describe the flow-subspace $S_{can}$, and end up with a regularity notion associated with a regular flow.

The pair $\{E, F\}$ is supposed to be \emph{pre-regular}
The first step of the reduction procedure from  \cite{RaRh} is then well-defined, we refer to \cite[Section 12]{RaRh} for the substantiating arguments.  Here we apply this procedure to homogeneous DAEs only.

We start by $E_0=E,\,F_0=F,\,m_0=m,\,r_{0}=r$, $\theta_0=\theta$, and consider the homogeneous DAE
\begin{align*}
 E_0x'+F_0 x=0.
\end{align*}
By means of a basis $Z_0:\mathcal I\rightarrow \Real^{m_0\times(m_0-r_0)}$ of $(\im E_0)^{\perp}=\ker E_0^{*}$ and a basis $Y_0:\mathcal I\rightarrow \Real^{m_0\times r_0}$ of $\im E_0$ we divide the DAE into the two parts
\begin{align*}
 Y_0^*E_0x'+Y_0^*F_0x=0,\quad Z_0^*F_0x=0.
\end{align*}
From  $\im[E_0,\,F_0]=\Real^m$ we derive that $\rank Z_0^*F= m_0-r_{0}$, and hence the subspace   
$S_{0}=\ker Z_0^*F$ has dimension $r_{0}$. Obviously, each solution of the homogeneous DAE must stay in the subspace $S_{0}$. Choosing a continuously differentiable basis 
 $C_0:\mathcal I\rightarrow \Real^{m_0\times r_0}$ of $S_{0}$, each solution of the DAE can be represented as $x=C_0 x_{[1]}$, with a function $x_{[1]}:\mathcal I\rightarrow \Real^{r_0}$ satisfying the reduced to size $m_1=r_{0}$ DAE
 \begin{align*}
 Y_0^*E_0C_0 x_{[1]}'+Y_0^*(F_0C_0+E_0C'_0)x_{[1]}=0.
\end{align*}
Denote $E_1=Y_0^*E_0C_0$ and $F_1=Y_0^*(F_0C_0+E_0C'_0)$ which have size $m_1\times m_1$.
The pre-regularity assures that $E_1$ has constant rank $r_{1}=r_0-\theta_0\leq r_{0}$. Namely,
we have
\begin{align*}
 \ker E_1=\ker E_0C_0=C^{+}_0 (\ker E_0\cap S_{0}),\quad \dim \ker E_1= \dim (\ker E_0\cap S_{0})=\theta_0.
\end{align*}

Next we repeat the reduction step supposing that the new pair is  pre-regular again, and so on.
This yields 
$ m\geq r_{0}\geq \cdots\geq r_{j}\geq r_{j-1}\geq\cdots \geq 0 $.
Denote by $\mu$ the smallest integer  such that either $r_{\mu-1}=r_{\mu}>0$ or $r_{\mu-1}=0$.
Then, it follows that $\ker E_{\mu-1}\cap S_{\mu-1}=\{0\}$,  which means in turn  that
\begin{align*}
 E_{\mu-1}x_{[\mu-1]}'+F_{\mu-1}x_{[\mu-1]}=0
\end{align*}
represents a regular index-1 DAE. If  $r_{\mu-1}=0$, that is $E_{\mu-1}=0$, then $F_{\mu-1}$ is nonsingular, which leads to a zero flow $x_{[\mu-1]}(t)\equiv 0$. On the other hand,
 if $r_{\mu}>0$ then $E_{\mu}$ remains nonsingular and 
\begin{align*}
 E_{\mu}x_{[\mu]}'+F_{\mu}x_{[\mu]}=0
\end{align*}
  is an implicit regular ODE living in $\Real^{m_{\mu}}$,\;$m_{\mu}=r_{\mu-1}$.
Finally, if $r_{\mu}>0$ each solutions of the original homogeneous DAE \eqref{DAE0} has the form
\begin{align*}
 x=\underbrace{C_0C_1\cdots C_{\mu-1}}_{=:C}x_{[\mu]}=:Cx_{[\mu]}, \quad C:[a,b]\rightarrow\Real^{m\times r_{\mu-1}},\; \rank C=r_{\mu-1}.
\end{align*}
Otherwise, if $r_{\mu}=0$, there is only the  identically vanishing solution of the homogeneous DAE, $x=0$.

Moreover, for each $\bar t\in \mathcal I$ and each $z\in\im C(\bar t)$, there is exactly one solution of the original homogeneous DAE passing through, $x(\bar t)=z$. 

As proved in \cite{RaRh}, the ranks $r=r_{0}> r_{1}>\cdots> r_{\mu-1}$ are independent of the special choice of the involved basis functions.

The property of pre-regularity does not necessarily carry over to the subsequent pair, as Example \ref{e.1} shows. 
\begin{example}\label{e.1}
The pair $\{E,F\}$ from \cite[p. 91]{GM86}, 
\begin{align*}
  E(t)=\begin{bmatrix}
        -t&t^2\\-1&t
       \end{bmatrix},\quad F(t)=\begin{bmatrix}
        1&0\\0&1
       \end{bmatrix},
\end{align*}
is pre-gegular with $m=2$, $r=1$ and $\theta =1$. Choosing bases 
\begin{align*}
 C_0=\begin{bmatrix}
      t\\1
     \end{bmatrix},\quad 
Y_0=\begin{bmatrix}
      t\\1
     \end{bmatrix}, \quad
Z_0=\begin{bmatrix}
       1\\-t
    \end{bmatrix},
\end{align*}
the subsequent pair is given by  $m_1=r=1$,  $E_1=Y^*_0EC_0=0$, and $F_1=Y^*_0FC_0+Y^*_0EC'_0=0$, such that $[E_1\,F_1]$ fails to have full row-rank $m_1$, thus the pair $\{E_1,F_1\}$ fails to be pre-regular. The associated to $\{E,F\}$ homogeneous DAE possess the solutions 
\begin{align*}
 x(t)= \alpha(t)\begin{bmatrix}
       t\\1
      \end{bmatrix},
\end{align*}
in which $\alpha$ stands for an arbitrary smooth real function, which does not fit our idea of regularity.
\qed
\end{example}

\begin{definition}\label{d.2}
 The pre-regular pair  $\{E,F\}$ with $r<m$  and the associated DAE \eqref{DAE0}, respectively, are called \emph{regular}\footnote{In \cite{RaRh} instead the term \emph{completely regular} is used whereas pre-regular pairs like example \ref{e.1} are called \emph{regular}.We have not taken up this notation, but that of other works, which seems more suitable to us.} if there is an integer $\mu\in\Natu$ such that 
 the above reduction procedure is well-defined up to level $\mu-1$, each pair $\{E_{i},F_{i}\}$, $i=0,\ldots,\mu-1$, is pre-regular,  and if $r_{\mu-1}>0$ then $E_{\mu}$ is well-defined and  nonsingular, $r_{\mu}=r_{\mu-1}$. If $r_{\mu-1}=0$ we set $r_{\mu}=r_{\mu-1}$.
 The index $\mu$ and the ranks $r=r_{0}> r_{1}>\cdots > r_{\mu-1}=r_{\mu}$ are called \emph{characteristic values} of the pair and the DAE, respectively.
\end{definition}
By construction, for a regular pair it follows that $r_{i+1}= r_{i}-\theta_{i}$, $i=0,\cdots,\mu-1$. Therefore, in place of the above $\mu+1$ rank values $r_0,\ldots,r_{\mu}$, the following rank and  dimensions,
\begin{align}\label{theta}
 r \quad \text{and}, \quad \theta_0 \geq\theta_1 \geq \cdots \geq \theta_{\mu-2} >\theta_{\mu-1}=0, 
\end{align}
can serve as characteristic quantities. Later it will become clear that these data also play a supporting role in other concepts, too.
\begin{remark}\label{r.RR}
A predecessor version of the reduction procedure in \cite{RaRh}
was already proposed and analyzed in \cite{Cis1982} under the name  \emph{elimination of the unknowns}, even for  more general pairs of rectangular matrix functions.
There, an additional scaling of the respective pairs is incorporated to put them in partitioned form on each stage, cf.\ Appendix~\ref{subs.CistJans}, but this makes the description less clear. The regularity notion given in \cite{Cis1982} is consistent with Definition \ref{d.2}. 
Another very related such reduction technique has been presented and extended a few years ago under the name \emph{dissection concept} \cite{Jansen2014}, see also Appendix~\ref{subs.CistJans}. This notion of regularity also agrees with Definition \ref{d.2}.
As we shall see below, the regularity notions related to the strangeness-index concept and the tractability-index framework  are consistent with Definition \ref{d.2}, too.
Furthermore, the understanding of regular points e.g. in \cite[Section 2.2.7]{RR2008} fits to this then also. 
\end{remark}

\begin{theorem}\label{t.Scan}
 If the DAE  \eqref{DAE0} is regular on $\mathcal I$ with index $\mu$ and  characteristics 
 \[r=r_{0}> r_{1}>\cdots> r_{\mu-1},\]
 or, equivalently, \eqref{theta},
  then $S_{can}(t)$ has  dimension $d=r-\sum_{j=0}^{\mu-2}\theta_j=r_{\mu-1}$ for all $t\in\mathcal I$, and
  the matrix function $C:\mathcal I\rightarrow \Real^{m\times d}$, $C=C_{0}\cdots C_{\mu-2}$, generated by the reduction procedure is a basis of $S_{can}$.
\end{theorem}
\begin{proof}
Regarding the relation $r_{i+1}= r_{i}-\theta_{i}$, $i=0,\cdots,\mu-2$ directly resulting from the reduction procedure, the assertion 
 is an immediate consequence of \cite[Theorem 13.3]{RaRh}.
\end{proof}
\begin{remark}\label{r.subspaces}
 The subspaces $S_{j}= \ker Z^*_jF_j=\im C_j\subset\Real^{m_j}$, with  $m_j=r_{j-1}=\rank E_{j-1}$,  $j=1,\ldots,\mu-1$,  are living in spaces of different dimension $m_j$. In contrast, letting 
 \begin{align*}
  S^{[0]}&:=S_0,\\
  S^{[j]}&:=\im C_0C_1\cdots C_j= C_0C_1\cdots C_{j-1}\im C_j  = C_0C_1\cdots C_{j-1} S_{j}, \; j=1,\ldots,\mu-1,
 \end{align*}
we arrive at the sequence of subspaces living all in $\Real^m$,
\begin{align*}
 S_0=S^{[0]}\supset S^{[1]}\supset\cdots\supset S^{[\mu-1]},
\end{align*}
showing dimensions $r=r_{0}>r_{1}>\cdots>r_{\mu-1}$, respectively.\footnote{Note that this sequence is
closely related to the first Wong chain in \cite{BeIlTr12}.}
Regarding that
\begin{align*}
 \ker E_0=\ker E,\quad
 \ker E_j=\ker EC_0\cdots C_{j-1}= (C_0\cdots C_{j-1})^+ \ker E,\; j=1,\ldots \mu-1,
\end{align*}
we derive the intersections $\ker E_j\cap S_{j}=(C_0\cdots C_{j-1})^+ \ker E\cap S^{[j]}$, 
and also the sequence of inclusions
\begin{align}\label{intersec}
  S_0\cap\ker E=S^{[0]}\cap\ker E\supset S^{[1]}\cap\ker E\supset\cdots\supset  S^{[\mu-1]}\cap\ker E,
\end{align}
with 
\begin{align*}
\dim S^{[j]}\cap\ker E&=\dim S_{j}\cap\ker E_j=\dim \ker E_{j+1}= r_{j}-r_{j+1}=\theta_j,\\j&=0,\ldots,\mu-1
\end{align*}
These inclusions seem to indicate a certain relationship to the  approach in \cite{EsLa2018}.
\qed
\end{remark}
It is important to mention that pre-regularity and regularity persist and the characteristic values are invariant under equivalence transformations. An equivalence transformation  of $\{E,F\}$ is given by matrix functions $L:\mathcal I\rightarrow \Real^{m\times m}$ being  continuous and pointwise nonsingular and $K:\mathcal I\rightarrow \Real^{m\times m}$ being continuously differentiable and pointwise nonsingular yielding the pair $\tilde E,\tilde F$,
\begin{align}\label{equivalence}
 \tilde E=LEK,\quad \tilde F=LFK+LEK'.
\end{align}
Obviously, if $\{E,F\}$ is pre-regular, then it follows $\tilde N=K^{-1}N$, $\tilde S=K^{-1}S$, and $\tilde N\cap \tilde S=K^{-1}(N\cap S)$, thus $\tilde r=r$, $\tilde{\theta}=\theta$. For the proof concerning regularity we refer to \cite{RaRh}.

\section{On further regularity notions and their relations}\label{s.comparison}
We are concerned here with the regularity notions and approaches from \cite{Cis1982,Jansen2014,KM2006,CRR} associated with the elimination procedure, the dissection concept, the strangeness reduction, and the tractability framework compared to Definition \ref{d.2}. 
The approaches in \cite{Cis1982,Jansen2014,KM2006} are 
de facto special solution methods including reduction steps by elimination of variables and differentiations of certain variables.
In contrast,  the concept in \cite{CRR} aims at a structural projector-based decomposition of the given DAE in order to analyze them subsequently.

We have already mentioned in Remark \ref{r.RR} that the elimination procedure in \cite{Cis1982} is an earlier, but less elegant version of the reduction procedure in \cite{RaRh}, which we have adopted in Section~\ref{s.regular}.
The regularity definition \cite[p. 58]{Cis1982} agrees with  Definition \ref{d.2} in the matter and also with the name. However, it does not yet specify any characteristic values, which is why we will only refer to \cite{RaRh} in the remainder of this paper.

Each of the concepts is associated with a sequence of pairs of matrix functions, each supported by certain rank conditions that look very different. 
Thus also the regularity notions, which require in each case that the sequences are well-defined with well-defined termination, look very different.
At the end of this section, we will know that all these regularity terms agree with our Definition \ref{d.2}, and that the characteristics \eqref{theta} capture all the rank conditions involved.

When describing the different methods, traditionally the same terms are used, for example $\{E_j, F_j \}$ for the matrix function pairs and  $r_j$ for the ranks. However, they have different meanings in each instance.
 To avoid confusion, we label the different characters with corresponding exponents $R$ (reduction), $S$ (strangeness), $D$ (dissection), and $T$ (tractability), respectively, when there is a risk of confusion. 
\bigskip

We first  relate  the regularity notion given by Definition \ref{d.2} to the  strangeness reduction concept. 

Within the strangeness reduction framework the following five rank-values of the matrix function pair $\{E,F\}$ play their role, e.g., \cite[p. 59]{KM2006}:
\begin{align}
 r&=\rank E,\label{S1}\\
 a&=\rank Z^*FT,\label{S2}\\
 s&=\rank V^*Z^*FT^{c},\label{S3}\\
 d&=r-s,\label{S4}\\
 u&=m-r-a-s,\label{S5}
\end{align}
whereby $T,T^{c}, Z,V$ represent orthonormal bases (ONBs) of $\ker E$, $(\ker E)^{\bot}$, $(\im E)^{\bot}$, and $(\im Z^*FT)^{\bot}$, respectively.
The strangeness concept is tied to the requirement that $r,a$, and $s$ are well-defined constant integers.
\begin{lemma}\label{l.prerKM}
 The pair $\{E,F\}$ is pre-regular, if and only if the values \eqref{S1}-\eqref{S5} are constant and $u=0$. In case of pre-regularity, one has
 \begin{align*}
  a=m-r-\theta,\quad s=\theta,\quad d=r-\theta.
 \end{align*}
\end{lemma}
\begin{proof}
 Let $\{E,F\}$ be pre-regular, and $T,T^c$ be smooth ONBs what is possible owing to the constant rank $r$ of $E$. Let $Z, V$ be pointwise ONBs and $N,S$ be given by \eqref{sub}. Then it holds $S=\ker Z^*F$ by construction, and $Z^*FT$ has size $(m-r)\times(m-r)$. From $\ker Z^*FT=T^*(N\cap S)$ we derive $\dim \ker Z^*FT=\theta$ and hence
 \begin{align*}
  a=\rank Z^*FT= m-r-\theta.
 \end{align*}
Next we investigate the strangeness value $s$. For this aim we decompose 
\[S=(N\cap S)\oplus X,
\]
with $X\subset S$, $\dim X=r-\theta$ and choose a basis $C_X$ of $X$ so that $\im C_X=X$, $\rank C_X=m-\theta$, $\im C_X\cap N=\{0\}$. 

Since $T^{c\,*}C_Xy=0$ means $C_Xy\in \ker E=N$, thus $C_Xy=0$, and therefore $y=0$, we learn that
\[\rank T^{c\,*}C_X=r-\theta.
\]
Regardind the decomposition
\[\Real^{m-r}=\im Z^*FT\oplus (\im Z^*FT)^{\perp}=\im Z^*FT\oplus \im V
\]
we conclude $\rank V=\theta$.

Now we inspect the nullspace of $V^* Z^*FT^c$. If $z\in \Real^r$ and 
$V^* Z^*FT^c z=0$, thus $Z^*FT^c z\in \im Z^*FT$,  then there is a $w\in \Real^{m-r}$ such that $Z^*FT^c z=- Z^*FTw$, further
\begin{align*}
 Z^*F[T^c\,T]\begin{bmatrix}
              z\\w
             \end{bmatrix}=0,\quad
             \text{that is}, \quad
             [T^c\,T]\begin{bmatrix}
              z\\w
             \end{bmatrix}\in S=(N\cap S)\oplus \im C_X.
\end{align*}
This yields $T^c z+Tw= \xi +C_X v$ with a $\xi\in N\cap S$ and a $ v\in \Real^{r-\theta}$. Taking into account that $T$ and $ T^c$ are ONBs of $N$ and $N^{\perp}$, respectively, we conclude $z=T^{c\,*}T^{c}z=T^{c\,*}C_Xv$, and hence
\[
\im T^{c\,*}C_X\subseteq \ker   V^* Z^*FT^c.
\]
On the other hand, choosing $z\in \im T^{c\,*}C_X$ we compute
\begin{align*}
 V^* Z^*FT^c z&=V^* Z^*FT^cT^{c\,*}C_X v\\
 &=V^* Z^*F(I-TT^*)C_X v\\
 &=V^* Z^*FC_X v-V^* Z^*FTT^*C_X v =0,
\end{align*}
yielding 
\begin{align*}
\im T^{c\,*}C_X= \ker   V^* Z^*FT^c,\quad \rank V^* Z^*FT^c=r-\rank T^{c\,*}C_X=r-(r-\theta)=\theta.
\end{align*}
Finally, we obtain $a+s=m-r-\theta+\theta=m-r$ and $u=0$. 
\medskip

On the other hand, let the pair $\{E,F\}$ have constant rank values \eqref{S1}--\eqref{S5}, and $u=0$.
Applying the  basic arguments of the strangeness reduction \cite[p. 68f]{KM2006} we transform the pair $\{E,F\}$  equivalently to $\{\mathring{\tilde E},\mathring{\tilde F}\}$,
 \begin{align*}
 \mathring{\tilde E}=\begin{bmatrix}
            I_s&&&\\&I_d&&\\
            &&0&\\
            &&&0
           \end{bmatrix},\quad   
 \mathring{\tilde F}=\begin{bmatrix}
            0&\tilde F_{12}&0&\tilde F_{14}\\
            0&0&0&\tilde F_{24}\\
            0&0&I_a&0\\
            I_s&0&0&0
           \end{bmatrix},        
 \end{align*}
with $d+s=r,\; a+s=m-r$.

For better clarity of the following we add a further permutation transformation. Using the permutation matrix
\begin{align*}
 K_P=\begin{bmatrix}
            0&0&0&I_s\\0&I_d&0&0\\
            0&0&I_a&0\\
            I_s&0&0&0
           \end{bmatrix}, 
\end{align*}
we arrive at
 \begin{align}\label{tildeEF}
 \tilde E=\mathring{\tilde E}K_P=\begin{bmatrix}
            0&0&0&I_s\\0&I_d&0&0\\
            0&0&0&0\\
            0&0&0&0
           \end{bmatrix},\quad   
  \tilde F=\mathring{\tilde F}K_P=\begin{bmatrix}
            \tilde F_{14}&\tilde F_{12}&0&0\\
            \tilde F_{24}&0&0&0\\
            0&0&I_a&0\\
            0&0&0&I_s
           \end{bmatrix},        
 \end{align}
which is, of course, again  equivalent to the original pair $\{E,F\}$. It remains checking its pre-regularity. The condition $\im [\tilde E \, \tilde F]=\Real^m$ is evident and $\tilde E$ has constant rank $r=s+d$. From 
\begin{align}\label{CP}
 \tilde N=\ker \tilde E=\im\begin{bmatrix}
                                I_s&0\\0&0\\0&I_a\\0&0
                               \end{bmatrix},\quad
\tilde S=\im \tilde C, \quad
                               \tilde C=\begin{bmatrix}
                                I_s&0\\0&I_d\\0&0\\0&0
                               \end{bmatrix},\quad
\tilde N\cap \tilde S=\im\begin{bmatrix}
                                I_s\\0\\0\\0
                               \end{bmatrix}.
\end{align}
we know that $\theta=s$, such that the pair $\{\tilde E,\tilde F\}$, and in turn the original pair $\{E,F\}$, is pre-regular.
\end{proof}
Assuming that the pair $\{ E, F\}$ is pre-regular, in turn the transformed pair $\{\tilde E,\tilde F\}$ given by \eqref{tildeEF} is also pre-regular, we  provide the successor pairs  $\{\tilde E^R_1,\tilde F^R_1\}$ according to the reduction procedure in Section \ref{s.regular} and $\{ E^S_1,F^S_1\}$ according to the strangeness framework \cite{KM2006}.

Using the basis functions $\tilde Y$ and $\tilde C$,
\begin{align*}
 \tilde Y^{*}=\begin{bmatrix}
               I_s&0&0&0\\0&I_d&0&0
              \end{bmatrix},\quad 
   \tilde C=\begin{bmatrix}
                                I_s&0\\0&I_d\\0&0\\0&0
                               \end{bmatrix},          
\end{align*}
we form the reduced pair $\{\tilde E^R_1,\tilde F^R_1\}$ to $\{\tilde E,\tilde F\}$in accordance with
Section \ref{s.regular}, that is,
\begin{align*}
 \tilde E^R_1=\tilde Y^{*}\tilde E\tilde C=\begin{bmatrix}
                                              0&0\\0&I_d
                                             \end{bmatrix},\quad
\tilde F^R_1=\tilde Y^{*}\tilde F\tilde C=\begin{bmatrix}
                                              \tilde F_{14}&\tilde F_{12}\\\tilde F_{24}&0
                                             \end{bmatrix}.
\end{align*}
The matrix function $\tilde E^R_1$ has size $r\times r$, $ r=s+d$  and constant rank $r_{[1]}=d$. If we look at the subspaces, 
\begin{align*}
 \tilde S^R_1=\{z\in \Real^{s+d}: \tilde F_{14}z_1+\tilde F_{12}z_2=0\},\quad \tilde S^R_1\cap \tilde E^R_1=\{z\in \Real^{s+d}: \tilde F_{14}z_1=0, z_2=0\},
\end{align*}
 it becomes clear that the reduced pair $\{\tilde E^R_1,\tilde F^R_1\}$ is pre-regular again, if and only if the following two conditions,
\begin{align}
 \im [\tilde F_{14}\, \tilde F_{12}]&=\Real^s, \label{B1}\\
 \theta_1:=\dim \tilde S^R\cap \ker \tilde E^R_1&=\dim \ker \tilde F_{14}\quad \text{is constant}.\label{B2}
\end{align}
are satisfied.
In particular, \eqref{B2} requires that the $s\times s$ matrix function $\tilde F_{14}$ shows constant rank $s-\theta_1=\theta-\theta_1$.
\bigskip

On the other hand, the matrix function pair $\{E^S_1, F^S_1\}$ following the original pair $\{E, F\}$  within  the strangeness framework is given by  replacing the 
entry $(\tilde E)_{14}$ in \eqref{tildeEF}  by a zero matrix, which leads to
\begin{align*}
  E^S_1=\begin{bmatrix}
            0&0&0&0\\0&I_d&0&0\\
            0&0&0&0\\
            0&0&0&0
           \end{bmatrix},\quad   
 F^S_1= \tilde F=\begin{bmatrix}
            \tilde F_{14}&\tilde F_{12}&0&0\\
            \tilde F_{24}&0&0&0\\
            0&0&I_a&0\\
            0&0&0&I_s
           \end{bmatrix}.        
 \end{align*}
We determine the strangeness characteristics of this pair. For this aim we use the corresponding bases
\begin{align*}
 T_1=\begin{bmatrix}
      I_s&0&0\\0&0&0\\0&I_a&0\\0&0&0
     \end{bmatrix},\quad
T^c_1=\begin{bmatrix}
      0\\I_d\\0\\0
     \end{bmatrix},\quad
Z_1=\begin{bmatrix}
      I_s&0&0\\0&0&0\\0&I_a&0\\0&0&I_s
     \end{bmatrix}.
\end{align*}
and obtain 
\begin{align*}
 Z^*_1F^S_1T_1=\begin{bmatrix}
              \tilde F_{14}&0&0\\0&I_a&0\\0&0&I_s
             \end{bmatrix},\quad
 Z^*_1F^S_1T^c_1=\begin{bmatrix}
              \tilde F_{12}\\0\\0
             \end{bmatrix},
\end{align*}
and further $r^S_1=d=r-\theta=r_{[1]}$ and $a_1= \rank Z^*_1F^S_1T_1=a+s+\rank \tilde F_{14}$. Set $r_{F14}=\rank F_{14}$.
Let $\tilde V_{14}$ denote an ONB of $(\im \tilde F_{14})^{\perp}$. $\tilde V_{14}$ has size $s\times (s-r_{F14})$. Then, the matrix function
\begin{align*}
 V_1=\begin{bmatrix}
      \tilde V_{14}\\0\\0
    \end{bmatrix},\quad \text{with size}\; (s+a+s)\times (s-r_{F14})
\end{align*}
forms an ONB of $(\im Z^*_1F^S_1T_1)^{\perp}$.
Next we compute 
\begin{align*}
 V^*_1Z^*_1F^S_1T^c_1=\tilde V^*_{14}\tilde F_{12},\quad s_1=\rank \tilde V^*_{14}\tilde F_{12}.
\end{align*}
further $u_1=m-r^S_1-a_1-s_1=m-r+s-(a+s+r_{F14})-s_1=s-r_{F14}-s_1$.
This results in the fact that the pair $\{E^S_1,F^S_1\}$ is pre-regular, 
exactly if $r_{F14}$ is constant, in turn $a_1$ is constant, and further $u_1=0$, so that $s_1=s-r_{F14}=\dim \ker \tilde F_{14}$. This requires exactly 
the conditions \eqref{B1} and \eqref{B2}. In the consequence, the pairs $\{\tilde E^R_1,\tilde F^R_1\}$ and $\{E^S_1,F^S_1\}$ are pre-regular simultaneously.
If they are pre-regular, then  $s_1=\theta_1,\;a_1=m-r +s-s_1=m-r+\theta-\theta_1$.
\begin{lemma}\label{l.3}
 Given are two matrix function pairs $\{E,F\}$ and $\{\bar E,\bar F\}$ acting in $\Real^m$ and   $\Real^{\bar m}$, respectively, and $\bar m<m$, 
 \begin{align*}
  E=\begin{bmatrix}
     \bar E&0\\0&0
    \end{bmatrix},\quad
F=\begin{bmatrix}
     \bar F&0\\0&I_{m-\bar m}
    \end{bmatrix}.
 \end{align*}
Then $\{E,F\}$ and $\{\bar E,\bar F\}$ are pre-regular simultaneously. If so, the characteristics are
\begin{align*}
 r=\bar r, \, a=\bar a+(m-\bar m),\; s=\bar s,\, \theta =\bar{\theta},\; s=\theta.
\end{align*}
\end{lemma}
\begin{proof}
Obviously, $E$ and $\bar E$ share their rank, and $\im [\bar E\; \bar F]=\Real^{\bar m}$  implies $\im [E\, F]=\Real^{m}$ and vice versa. From
\begin{align*}
\bar S&=\{w\in \Real^{\bar m}: \bar Fw\in \im \bar E\},\\
 S&=\{z\in \Real^{\bar m+(m-\bar m)}: \bar Fz_1\in \im \bar E, z_2=0\}=\bar S\times \{0\},\\
 N&=\bar N\times\Real^{m-\bar m},\\
 N\cap S&=(\bar N\cap \bar S)\times  \{0\},
\end{align*}
we conclude that the pairs are pre-regular simultaneously. Is they are so, then obviously, $r=\bar r$, $\theta=\bar \theta$.

If $\bar T,\, \bar T^c ,\, \bar Z$ are ONBs to $\ker \bar E,\, (\ker \bar E)^{\perp},\, (\im\bar E)^{\perp}$ then 
\begin{align*}
 T=\begin{bmatrix}
    \bar T&0\\
    0&I_{m-\bar m}
   \end{bmatrix},\quad 
   T^c=\begin{bmatrix}
    \bar T^c\\
    0
   \end{bmatrix},\quad
 Z=  \begin{bmatrix}
    \bar Z&0\\
    0&I_{m-\bar m}
   \end{bmatrix},\quad
\end{align*}
are ONBs to $\ker E,\, (\ker E)^{\perp},\, (\im E)^{\perp}$. This leads to
\begin{align*}
 Z^*FT=\begin{bmatrix}
    \bar Z^*\bar F\bar T&0\\
    0&I_{m-\bar m}
   \end{bmatrix},\quad a=\rank Z^*FT= \rank \bar Z^*\bar F\bar T+m-\bar m=\bar a+m-\bar m.
\end{align*}
Finally, owing to Lemma \ref{l.prerKM} we find $s=\theta$, $\bar s=\bar\theta$, and hence $s=\bar s$.
\end{proof}
\begin{theorem}\label{t.indexrelation}
Let the  pair $\{E,F\}$ be regular on $\mathcal I$ with index $\mu\in \Natu$ and  characteristics $r<m$,
$\theta_0=0$ if $\mu=1$, and, for $\mu>1$,
\begin{align*}
r<m,\quad \theta_0\geq\cdots\geq\theta_{\mu-2}>\theta_{\mu-1}=0.
\end{align*}
Then the following statements are true: 
\begin{description}
 \item[\textrm{(a)}] The strangeness index $\mu_S$ is well-defined for $\{E,F\}$, and $ \mu_S=\mu-1$. The associated characteristics are\footnote{The quantities $a_i, d_i, s_i, u_i$ only occur in the strangeness concept and do not require the $S$ label.}
 \begin{align*}
  r^S_0&=r,\\
  s_i&=\theta_i,\\
d_i&=r^S_i-\theta_i,\\
  a_i&=m-r_i-\theta_i,\\
  u_i&=0\\
  r^S_{i+1}&=d_i,\quad i=0,\ldots,\mu-1.
 \end{align*}
\item[\textrm{(b)}] The pair $\{E,F\}$ is regular with tractability index $\mu$ and characteristics 
\begin{align*}
 r^T_0=r,\quad r^T_i=m-\theta_{i-1},\quad i=1,\ldots,\mu.
\end{align*}
\item[\textrm{(c)}] The pair $\{E,F\}$ is regular with dissection index  index $\mu$ and characteristics 
\begin{align*}
 r^D_0=r,\quad r^D_i=m-\theta_{i-1},\quad i=1,\ldots,\mu.
\end{align*}
\end{description}
\end{theorem}
\begin{proof}
 We emphasize again that all characteristics involved here do not change under equivalence transformations \cite{RaRh,KM2006,CRR,Jansen2014}. 
 
 (a): We perform the reduction process of the strangeness concept step by step and compare each level with the reduction from \cite{RaRh}. With the help of Lemmata \ref{l.prerKM} and \ref{l.3} we get step by step the assertion. The reduction pairs in \cite{RaRh}, which have lower dimension, turn out to be (to equivalence exactly) those parts of the pairs from \cite{KM2006} which still play a role for the further.
 
 (b): Given Assertion (a), Assertion (b) is a direct consequence of \cite[Theorem 2.79]{CRR}. We refer to Appendix~\ref{AdmissibleMatrix} for  informations concerning the admissible matrix function sequences associated with the tractability index.
 
 (c): Given Assertion (b), Assertion (c) is a direct consequence of \cite[Theorem 4.25]{Jansen2014}. We refer to Appendix~\ref{subs.CistJans} for a brief description of the basic step in the dissection index concept. 
\end{proof}
\begin{corollary}\label{c.pencil}
 Let the pair $\{E,F\}$ of matrices $E,F\in \Real^{m\times m}$ be regular in the sense of Definition \ref{d.2} with characteristic values \eqref{theta}. Then the matrix pencil $\lambda E+F$ is regular and, for $i=0,\ldots,\mu-2$, the quantity $\theta_i$ is the number of Jordan blocks of order $\geq 2+i$ within the nilpotent  matrix $\mathcal N$ in the Weierstra{\ss}--Kronecker form \eqref{MP1}. 
\end{corollary}
\begin{proof}
 A constant matrix pair and the associated matrix pencil are known to be regular with Kronecker index $\mu$ exactly if it is regular with tractability index $\mu$, e.g.,\cite[Chapter 1]{CRR}, and if so, $l_{i}= m-r^T_{i-1}$ is the number of Jordan blocks of order  $\geq i$ in the nilpotent matrix $\mathcal N$ in \eqref{MP1}, $i=0,\ldots,\mu-2$. Owing to Theorem \ref{t.indexrelation} we find that  $l_i=m-r^T_{i-1}=m-(m-\theta_{i-2})=\theta_{i-2}$, for $i=2,\ldots,\mu$.
 \end{proof}
\begin{remark}\label{r.theta}
 It seems to us very worth highlighting, that the $\theta$-characteristics \eqref{theta} make sense in all approaches and allow to determine all further characteristic values. Additionally, Corollary \ref{c.pencil} reveals a feature independent of any method.
 Furthermore, so far it is clear that $\theta_i$ coincides with $s_i$ and it is the dimension of intersecting subspaces in different approaches\footnote{See Appendix~\ref{AdmissibleMatrix} for the definition of $G_i, B_i$.}:
 \begin{align*}
  &\ker E^R_i\cap S^R_i=\ker E^R_i\cap \{z\in \Real^{m_i}:F^R_i z\in \im E^R_i\},\\
  &\ker E\cap S^{[i]},\quad \text{in Remark \ref{r.subspaces}},\\
  &\ker G_i\cap S^T_i=\ker G_i\cap \{z\in \Real^{m}:B_i z\in \im G_i\}=\ker G_i\cap \{z\in \Real^{m}:B_0 z\in \im G_i\}.
 \end{align*}
\end{remark}
\begin{remark}\label{r.compare}
 We have compared here the application of the reduction procedures from \cite{KM2006,Jansen2014} and \cite{RaRh} to homogeneous DAEs, only. If, on the other hand, 
 inhomogeneities $q$ shall be considered, differentiations with respect to $q$ must be made at each level. 
In \cite{RaRh}, the resulting explicit relations are considered as being finished and, thus, neglected. So only the system of lower dimension, which is still of interest, is treated further. In contrast, in \cite{KM2006}, all equations are further carried along. 
 Apart from equivalence transformations at each level, this makes all the difference.
\end{remark}
\section{Accurate initial condition and the second canonical subspace $N_{can}$}\label{s.accurateIC}
We emphasize again that we are not looking for consistent initial values here, but for an adequate formulation of initial conditions that lead to uniquely solvable IVPs, i.e., that precisely determine the free parameters of the flow of the DAE.
We adopt the notion of  accurately stated boundary condition \cite[Definition 2.3]{LMW} for this purpose.
Consider the IVP
\begin{align}\label{IVP0}
 Ex'+Fx=q,\quad G_a x(a)=\gamma,
\end{align}
for a regular DAE featuring dynamical degree of freedom $d$, with  a matrix $G_{a}\in\Real^{s\times m}$, $s\geq d$, $a\in \mathcal I$, $\gamma\in \im G_a\subseteq\Real^s$.

Let $x_*$ be a solution of the IVP \ref{IVP0}. Following \cite[Definition 2.3]{LMW}, the initial condition in \ref{IVP} is \emph{accurately stated}, if all slightly perturbed IVPs
\begin{align}\label{IVP1}
 Ex'+Fx=q,\quad G_a x(a)=\gamma+\Delta\gamma,\quad \Delta\gamma\in \im G_a,
\end{align}
are uniquely solvable and their solutions satisfy, on a compact interval $\mathcal I_a\subseteq\mathcal I$, the inequality
\begin{align}\label{IVP2}
 \max_{t\in\mathcal I_a}\lvert x(t)-x_*(t)\rvert \leq K \lvert \Delta\gamma \rvert .
\end{align}
It should be remembered that regular  higher-index DAEs lead to ill-posed problems, even if the initial conditions are stated accurately, e.g.,\cite{Mae2014}. 
\bigskip

Regular index-1 DAEs are studied in detail in \cite{GM86,CRR}, cf. also Subsection \ref{s.semi}.
Their dynamical degree of freedom is $d=r=\rank E(t)$, and  one has simply
\[N_{can}=N,\quad S_{can}=S
\]
with $N$ and $S$ from \eqref{sub}.
The related canonical projector function $\pPi_{can}:\mathcal I\rightarrow \Real^{m\times m}$ is given and also the solvability statements for IVPs with the initial condition
\begin{align*}
 x(a)-x_{a}\in  N(a)=N_{can}(a), \quad  x_a\in \Real^m,
\end{align*}
 are proved. We underline that here $x_a$ is arbitrary and it is not necessarily a consistent value. On the other hand, $\pPi_{can}(a)x(a)=\pPi_{can}(a)x_a$ is always valid for the solution, while $x(a)=x_a$ cannot be expected in general. 

While in the index-1 case one has simply $N_{can}=N$ and $S_{can}=S$, in the case of regular higher-index pairs   the canonical subspaces will be subspaces with higher and lower dimension, respectively, and
\begin{align*}
 N_{can}\supset N,\quad S_{can} \subset S.
\end{align*}
\begin{theorem}\label{t.Ncan}
 Let the  pair $\{E,F\}$ be regular on $\mathcal I$ with index $\mu\in \Natu$ and  characteristics $r<m$,
 $\theta_0=0$ if $\mu=1$, and, for $\mu>1$,
\begin{align*}
 r<m, \quad \theta_0\geq\cdots\geq\theta_{\mu-2}>\theta_{\mu-1}=0.
\end{align*}
Let $S_{can}$ be the canonical flow subspace provided by Theorem \ref{t.Scan}. Then the following statements are true: 
\begin{description}
\item[\textrm{(a)}] There is a  subspace $N_{can}\subset\Real^m$  having, on the given Interval $\mathcal I$, the dimension $m-d=m-r+\sum_{j=0}^{\mu-2}\theta_j$, and the decomposition
\begin{align*}
 \Real^m=S_{can}\oplus N_{can}
\end{align*}
holds pointwise on all $\mathcal I$. Additionally, the canonical projector $\pPi_{can}$ onto $S_{can}$ along $N_{can}$ is well-defined on all $\mathcal I$.
\item[\textrm{(b)}]  If $G_a\in \Real^{s\times m}$ is such that $s\geq d$, $\ker G_a=N_{can}$, and if $q$ is sufficiently smooth, then each IVP
\begin{align*}
 Ex'+Fx=q,\quad G_{a}(x(a)-x_a)=0, \quad x_a\in \Real^m,\quad a\in\mathcal I,
\end{align*}
has a unique solution which features the relation  
$\pPi_{can}(a)x(a)=\pPi_{can}(a)x_a$ and satisfies, on compact intervals $\mathcal I_a\subseteq\mathcal I$ around $a$, the inequality
\begin{align}\label{inequality}
 \max_{t\in\mathcal I_a}\lvert x(t)\rvert \leq K\{ \lvert \pPi_{can}(a)x_a \rvert + \max_{t\in\mathcal I_a}\lvert q(t)\rvert +\sum^{\mu-1}_{l=1} \max_{t\in\mathcal I_a}\lvert q^{(l)}(t)\rvert\},
\end{align}
with a constant $ K$ depending on the pair $\{E,F\}$ and the interval $\mathcal I_a$ only.
\end{description}
\end{theorem}
The assertions justify the notation $N_{can}$.
\begin{proof}
Owing to Theorem \ref{t.indexrelation} the pair  $\{E,F\}$ is regular with tractability index $\mu$ and characteristics 
\begin{align*}
 r^T_0=r,\quad r^T_i=m-\theta_{i-1},\quad i=1,\ldots,\mu.
\end{align*}
This means in the projector-based framework that there is an admissible matrix function sequence (see \cite{CRR}, also Appendix~\ref{AdmissibleMatrix} below). The related nullspaces $N_i$ have dimensions $\theta_{i-1}$, $i=1,\cdots,\mu$, and $\dim N_0=m-r$, and
\begin{align}\label{N1}
 \ker \pPi_{\mu-1}=N_0+N_1+\cdots+N_{\mu-1}=N_0\oplus N_1\oplus\cdots\oplus N_{\mu-1}=:N_{can}
\end{align}
has dimension $\dim N_{can}=\sum_{i=0}^{\mu-1} \dim N_i =m-r+ \sum_{i=1}^{\mu-1} \dim \theta_{i-1}$. The subspace $N_{can}$ is shown to be independent of the special choice of the admissible projector functions which form the projector function $\pPi_{\mu-1}$, \cite[Theorem 2.8]{CRR}. Since we suppose sufficiently smooth $E$ and $F$, a so-called fine decoupling sequence can be constructed starting with an arbitrary projector function  $\pPi_0$ onto $N_0=\ker E$. Then, for $\pPi_{\mu-1}$ associated with a fine decoupling sequence, there is a further special projector function 
$Q_{* 0}$ onto $N_0$, such that (cf. \cite[Lemma 2.31, Theorem 2.42]{CRR}) $\pPi_{can}:= (I-Q_{* 0})\pPi_{\mu-1}$ is again a projector function. $\pPi_{can}$ has constant rank $d$ on all $\mathcal I$, and by construction $\ker\pPi_{can}=N_{can}$. It is further verified in \cite[Section 2.6]{CRR} that  $\im\pPi_{can}=S_{can}$, which in turn shows that $N_{can}$ is actually a complementary subspace to $ S_{can}$.

Assertion (b) is now completely verified by \cite[Theorem 2.63]{CRR}.
\end{proof}
\begin{corollary}\label{c.1}
 Under the conditions of Theorem \ref{t.Ncan} concerning the pair  $\{E,F\}$ and the matrix $G_a$, the initial condition in the IVP \eqref{IVP0} is accurately stated.
\end{corollary}
\begin{proof}
 Since $\gamma, \Delta\gamma\in\im G_a$, we may choose $x_a, \Delta x_a\in \Real^m$ such that $\gamma=G_a x_a$, $\Delta\gamma=G_a \Delta x_a$, and the IVPs \eqref{IVP1} are uniquely solvable by Theorem \ref{t.Ncan}(b). Furthermore, $z:=x-x_*$ satisfies the homogeneous DAE and the initial condition $G_az(a)-G_a\Delta x_a=0$, and \eqref{inequality} implies 
 \begin{align*}
   \max_{t\in\mathcal I_a}\lvert z(t)\rvert \leq K\ \lvert \pPi_{can}(a)\Delta x_a \rvert=K\ \lvert G_a^-G_a\Delta x_a \rvert\leq \tilde K \lvert G_a\Delta x_a \rvert= \tilde K \lvert \Delta \gamma \rvert,
 \end{align*}
what was to show.
\end{proof}

Finally in this part, let us add that the maximal fundamental solution matrix of a  regular DAE normalized at $a\in \mathcal I$, that is the solution of the IVP,
\begin{align*}
 EX'+FX=0,\quad X(a)=\pPi_{can}(a)
\end{align*}
feature semigroup properties and $\im X(t,a)=S_{can}(t)$, $\ker  X(t,a)=N_{can}(a)$, what we could see so similarly also in the case of constant matrix pairs in Section \ref{s.pencil}.

\section{An useful representations of the matrix $G_a$ for accurately stated initial conditions and the projector function \protect$\pPi_{can}\protect$}\label{subs.ker}
In the framework of the projector-based analysis \cite{CRR} admissible matrix function sequences and incorporated admissible projector functions play their role, see also Appendix~\ref{AdmissibleMatrix}.
Owing to \cite[Theorem 2.8]{CRR} the subspaces $N_{0}+ N_{1}+\cdots+ N_{i}$ do not depend of the special choice of the involved projector functions.
For a  regular DAE with  index $\mu$, it holds that
\begin{align*}
 \ker \pPi_{can}= \ker \pPi_{\mu-1}=N_{0}\oplus N_{1}\oplus\cdots\oplus N_{\mu-1},
\end{align*}
with $\pPi_{\mu-1}$ given by an arbitrary admissible matrix function sequence.  However, the possibilities of practical calculation are still limited, so we are looking for another way.
\medskip

The homogeneous adjoint DAEs to the above DAEs \eqref{sDAE} and \eqref{pDAE}, that is, $Ex'+Fx=q$ and $A(Dx)'+Bx=q$, are 
\begin{align}
-E^*y'+(F^*-{E^{*}}')y=0, \label{sDAEadj}\\
 -D^*(A^*y)'+B^*y=0. \label{pDAEadj}
\end{align}
Owing to \cite[Theorem 3 and Corollary 2]{LinhMae}, the original DAE and its adjoint are regular with  index $\mu$ at the same time, and they share the related characteristics, in particular the dynamical degree of freedom $d$.

Let $\pPi_{*\,can}$ denote the canonical projector function associated with the adjoint DAE and let $C_{*}:\mathcal I\rightarrow\Real^{m\times d}$ be a basis of $S_{*\,can}=\im \pPi_{*\,can}$, which can be provided, for example, using the procedure from Section \ref{s.regular}.

Let's take a brief look at
the semi-explisit index-1 DAE in Subsection \ref{s.semi} and consider its adjoint DAE.
\begin{example}\label{e.1cont}
The pair $\{E_*,F_*\}$ below describes the DAE adjoint to the DAE in \ref{s.semi},
 \begin{align*}
 E_*=-E^*=-E=-\begin{bmatrix}
    I_r& 0\\0&0
   \end{bmatrix},\quad
F_*=F^*= \begin{bmatrix}
    F^*_{11}& F^*_{21}\\F^*_{12}&F^*_{22}
   \end{bmatrix},\quad \text{with } F_{22}\; \text{ remaining  nonsingular}.
 \end{align*}
The pair is regular with index one, 
with \begin{align*}
      N_*=\im \begin{bmatrix}
             0\\I_{m-r}
            \end{bmatrix},\quad
S_*=\im C_*,\quad C_*= \begin{bmatrix}
             I_r\\-F_{22}^{*\,-1}F^*_{12}
            \end{bmatrix}, \quad N_*\cap S_*=\{0\},
     \end{align*}
such that $\Real^m=S_*\oplus N_*$. The canonical projector function $\pPi_{*\,can}$ onto $S_*$ along $N_*$ and the orthogonal projector $P_*$ onto $N_{*}^{\perp}$ along $N_*$  are given by
\begin{align*}
 \pPi_{*\,can}=\begin{bmatrix}
             I_r&0\\-F_{22}^{*\,-1}F^*_{12}&0
            \end{bmatrix},\quad
P_*=\begin{bmatrix}
             I_r&0\\0&0
            \end{bmatrix}.
\end{align*}
Observe that evidently here $d=r$ and 
\begin{align*}
C^*_{*}E&= ( I_d\, -F_{12}F_{22}^{-1}) E= [I_d\, 0],\\
\ker \pPi_{can}&=\ker C^*_{*}E,
\end{align*}
and hence, for stating accurate initial conditions one can choose 
\begin{align*}
 G_a (x(a)-x_a)=0,\quad G_a=C^*_{*}(a)E(a)
\end{align*}
However, in the index-1 case in contrast to all higher-index cases  there is the simpler possibility for the given DAE, namely, by choosing $x(a)-x_a\in \ker E=N=N_{can}$.
\qed
\end{example}

\begin{theorem}\label{t.ker}
 Let the pair $\{E, F\}$ be regular with index $\mu$ and canonical subspaces $S_{can}$ and  $N_{can}$, with bases
 $C_{S_{can}}$ and $C_{N_{can}}$, respectively.
 
 Then the adjoint pair  $\{-E^*, F^*-{E^*}'\}$ is also regular with the same characteristics. Moreover, with bases $C_{S_{*\,can}}$, $C_{N_{*\,can}}$ of their canonical 
 subspaces $S_{*\,can}$ and  $N_{*\,can}$, it results  that 
 \begin{align*}
 N_{can}=\ker C_{S_{*\,can}}^* E, \quad N_{*\,can}=\ker C_{S_{can}}^* E^*.
\end{align*}
\end{theorem}
\begin{proof}
We emphasize again that with the smoothness assumed here in general, it does not matter whether one chooses a homogeneous DAE in standard form, $Ex'+Fx=0$, or a DAE with proper involved derivative, $A(Dx)'+Bx=0$, to a given regular pair $\{E,F\}$, and $E=AD$, $B=F_AD'$, see Appendix~\ref{subs.forms}.
The pleasant  symmetry of the DAE with proper involved derivative to its adjoint in many cases facilitates the investigation, although $-D^*(A^*y)'+B^*y=0$ is only another notation for $-E^*y'+(F^*-{E^*}')y=0$.

By \cite[Theorem 3]{LinhMae}, the adjoint pair inherits the regularity along with characteristics from the regular pair $\{E,F\}$.

Owing to \cite[Lemma 3]{LinhMae} we know that
\begin{align*}
 D\pPi_{can}D^-&=(A^*\pPi_{*\,can}A^{*\,-})^*, 
\end{align*}
in which  $D^-$ and $A^{*\,-}$ are special generalized inverses corresponding to so-called complete decouplings. Regarding the relations 
\begin{align*}
 \ker D\pPi_{can}=\ker \pPi_{can},\quad \im\pPi_{* \, can} A^{*\,-}=\im\pPi_{* \, can},
\end{align*}
we derive
\begin{align*}
 \ker \pPi_{can}= \ker (A^{*\,-})^*\pPi_{*\,can}^* AD=\ker \pPi_{*\,can}^* AD.
\end{align*}
On the other hand, writing shorter $C_{*}=C_{S_{*\,can}}$  it holds that
\begin{align*}
 \ker C_*^*=(\im C_*)^{\perp}=(\im \pPi_{*\,can})^{\perp}=\ker \pPi_{*\,can}^*,
\end{align*}
further
\begin{align*}
 \ker C_*^*AD&=(\im (AD)^*C_*)^{\perp}=((AD)^*\im \pPi_{*\,can})^{\perp}(\im (AD)^*\pPi_{*\,can})^{\perp}\\
 &=\ker \pPi_{*\,can}^*AD,
\end{align*}
and hence $\ker \pPi_{can}=\ker C_*^*AD$. The second relation in the assertion of the theorem is valid for symmetry arguments.
\end{proof}

Having the basis   $C=:C_{S_{can}}$ of $S_{can}$, and a basis $C_{N_{can}}$ of  $\ker C_*^*AD=N_{can}$, respectively, the $m\times m$ matrix function $\mathcal M:=[C_{S_{can}}\,C_{N_{can}}]$  remains nonsingular everywhere on $\mathcal I$, and the canonical projector is given by
\begin{align*}
 \pPi_{can}=\mathcal M
 \begin{bmatrix}
             I_d&\\&0_{m-d}
            \end{bmatrix} \mathcal M^{-1}.
\end{align*}

\section{As an illustrating example: The linearized Campbell--Moore DAE}\label{s.CamobellMoore}
We investigate the linear DAE $E(t)x'(t)+F(t)x(t)=q(t)$ given by its coefficients
\begin{align*}
E&=\begin{bmatrix}
    1&0&0&0&0&0&0\\
    0&1&0&0&0&0&0\\
    0&0&1&0&0&0&0\\
    0&0&0&1&0&0&0\\
    0&0&0&0&1&0&0\\
    0&0&0&0&0&1&0\\
    0&0&0&0&0&0&0
   \end{bmatrix}=
   \begin{bmatrix}
    I_{3}&0&0\\
    0&I_{3}&0\\
    0&0&0
   \end{bmatrix},\\   
 F(t)&=\begin{bmatrix}
       0&0&0&-1&0&0&0\\
       0&0&0&0&-1&0&0\\
       0&0&0&0&0&-1&0\\0&0& \sin t&0&1&- \cos t&-2\rho \cos^{2} t\\
       0&0&-\cos t&-1&0&-\sin t&-2\rho\sin t\cos t\\
       0&0&1&0&0&0&2\rho\sin t\\
       2\rho\cos^2 t&2\rho\sin t\cos t&-2\rho\sin t&0&0&0&0
      \end{bmatrix}\\
      &=
      \begin{bmatrix}
       0&-I_{3}&0\\
       \mathfrak{C}(t)&\mathfrak{A}(t)&\mathfrak{B}(t)\\
       -\mathfrak{B}(t)^*&0&0
      \end{bmatrix}, \quad \rho\neq 0,
\end{align*}
\begin{align*}
\mathfrak{A}(t)= \begin{bmatrix}
                0&1&-\cos t\\
                -1&0&-\sin t\\
                0&0&0
               \end{bmatrix},\quad
  \mathfrak{B}(t)= 2\rho \begin{bmatrix}
                -\cos^2 t\\
                -\sin t\cos t\\
                \sin t
               \end{bmatrix},\quad
\mathfrak{C}(t)= \begin{bmatrix}
                0&0&\sin t\\
                0&0&-\cos t\\
                0&0&1
               \end{bmatrix},
\end{align*}
whereby the resulting blocks of different sizes with only zero entries are all denoted by $0$. 
Emphasizing the block structure we may  write the homogeneous DAE as
\begin{align}\label{DAE1}
 \begin{bmatrix}
  X'_{1}\\
   X'_{2}\\
   0
 \end{bmatrix}+
 \begin{bmatrix}
       0&-I_{3}&0\\
       \mathfrak{C}&\mathfrak{A}&\mathfrak{B}\\
       -\mathfrak{B}^*&0&0
      \end{bmatrix} \begin{bmatrix}
  X_{1}\\
   X_{2}\\
   X_{3}
 \end{bmatrix}= 0, \quad X_{1}=\begin{bmatrix}
  x_{1}\\
   x_{2}\\
   x_{3}
 \end{bmatrix},\; X_{2}=\begin{bmatrix}
  x_{4}\\
   x_{5}\\
   x_{6}
 \end{bmatrix},\; X_{3}=x_{7}.
\end{align}
Swapping the first two lines and also the variables leads to
\begin{align*}
 \begin{bmatrix}
  X'_{2}\\
   X'_{1}\\
   0
 \end{bmatrix}+
 \begin{bmatrix}
  \mathfrak{A}&\mathfrak{C}&\mathfrak{B}\\
       -I_{3}&0&0\\
       0&-\mathfrak{B}^*&0
      \end{bmatrix} \begin{bmatrix}
  X_{2}\\
   X_{1}\\
   X_{3}
 \end{bmatrix}= 0, \quad \text{with}\quad \mathfrak{B}^*\mathfrak{B}= 4\rho^2>0,
\end{align*}
what shows Hessenberg structure of size three. The DAE \eqref{DAE1} is regular with index $\mu=3$ and features the sizes $m=7$, $ m_{1}=3,\, m_{2}=3,\, m_{3}=1$,  dynamical degree of freedom $d=4$, as well as characteristical values $ r^T_{0}=r^T_{1}=r^T_{2}=6,\,r^T_{3}=7$.\footnote{See \cite[Theorem 3.42]{CRR}} In terms of Definition \ref{d.2} and \eqref{theta} one has $r=6$, $\theta_0=\theta_1=1$, and $\theta_2=0$.

The orthoprojector function $\Omega=\mathfrak{B}(\mathfrak{B}^*\mathfrak{B})^{-1}\mathfrak{B}^*=\frac{1}{4\rho^2}\mathfrak{B}\mathfrak{B}^* $ projecting pointwise the space  $\Real^3$ onto $\im \mathfrak{B}\subset \Real^3$ plays its role below. We have 
$\ker \Omega=\ker \mathfrak{B}^*=(\im \mathfrak{B})^{\bot}$.  

We choose a smooth basis $C_{\mathfrak{B}}$ of $\ker\mathfrak{B}^*$ such that
\begin{align}
 &\im C_{\mathfrak{B}}=\ker \mathfrak{B}^*=(\im \mathfrak{B})^{\bot}=\im(I_3-\Omega)\subset\Real^3,\quad \rank C_{\mathfrak{B}}=2, \label{xx} \\
 &\im \mathfrak{B}=\im \Omega,\quad\rank \Omega=1. \nonumber
\end{align}
It is evident that
\begin{align*}
 S=\im\begin{bmatrix}
       C_{\mathfrak B}&0&0\\
       0&I_3&0\\
       0&0&1
      \end{bmatrix},\quad 
N=\im\begin{bmatrix}
       0\\
       \vdots\\
       1
      \end{bmatrix},\quad S\cap N=N, \quad \dim S\cap N=1.
\end{align*}
By calculating the general solution of the DAE \eqref{DAE1} and inspecting the solution structure one obtains that
\begin{align*}
 S_{can}=\im C_{S_{can}}\subset S,\quad C_{S_{can}}=\begin{bmatrix}
                                            C_{\mathfrak B}&0\\
                                            -\Omega'C_{\mathfrak B}&C_{\mathfrak B}\\
                                            \mathfrak K_{1}&\mathfrak K_{2}
                                           \end{bmatrix},
\end{align*}
with $\mathfrak K_{1}=\frac{1}{4\rho^2}\mathfrak B^*(2\Omega'-\mathfrak A)C_{\mathfrak B}$, \, $\mathfrak K_{2}=\frac{1}{4\rho^2}\mathfrak B^*(\Omega''-\mathfrak C+\mathfrak A \Omega'-2\Omega'\Omega')C_{\mathfrak B}$.

\subsection{An admissible matrix function sequence for DAE \eqref{DAE1}}\label{subs:sequence}
We start by
\begin{align*}
 \pPi_{0}=P_{0}=G_{0}=AD=\begin{bmatrix}
       I_{3}&0&0\\
       0&I_{3}&0\\
       0&0&0
      \end{bmatrix},\quad B_{0}=\begin{bmatrix}
       0&-I_{3}&0\\
       \mathfrak{C}&\mathfrak{A}&\mathfrak{B}\\
       -\mathfrak{B}^*&0&0
      \end{bmatrix},\quad Q_{0}=I_{7}-P_{0},
\end{align*}
yielding
\begin{align*}
 G_{1}=G_{0}+B_{0}Q_{0}=\begin{bmatrix}
       I_{3}&0&0\\
       0&I_{3}&\mathfrak{B}\\
       0&0&0
      \end{bmatrix},\quad Q_{1}=\begin{bmatrix}
       0&0&0\\
       0&\Omega&0\\
       0&-\frac{1}{4\rho^2}\mathfrak{B}^*&0
      \end{bmatrix},
\end{align*}
and further
\begin{gather*}
 \pPi_{0}Q_{1}=\begin{bmatrix}
       0&0&0\\
       0&\Omega&0\\
       0&0&0
      \end{bmatrix},\quad \pPi_{1}=\begin{bmatrix}
       I_{3}&0&0\\
       0&I_3-\Omega&0\\
       0&0&0
      \end{bmatrix},\quad D \pPi_{1}D^{+}=\begin{bmatrix}
       I_{3}&0\\
       0&I_3-\Omega\\
      \end{bmatrix},
\\
\begin{aligned}
 G_{2}&=G_{1}+B_{0}\pPi_{0}Q_{1}-G_{1}D^{+}(D \pPi_{1}D^{+})'D\pPi_{0}Q_{1}\\
 &=
 \begin{bmatrix}
       I_{3}&-\Omega&0\\
       0&I_{3}+\mathfrak{A}\Omega+\Omega'\Omega&\mathfrak{B}\\
       0&0&0
      \end{bmatrix}.
\end{aligned}
\end{gather*}
Next we derive
\begin{align*}
Q_{2}&=\begin{bmatrix}    
       \Omega&0&0\\
      (I_3-(I_3-\Omega)(\mathfrak{A}+\Omega'))\Omega&0&0\\
      \mathfrak{D}&0&0
      \end{bmatrix},
\end{align*}
with $\mathfrak{D}=-\frac{1}{4\rho^2}\mathfrak{B}^{*}(I_3+\mathfrak{A}\Omega+\Omega'\Omega)
(I_3-(I_3-\Omega)\mathfrak{A}\Omega-\Omega'\Omega)
=-\frac{1}{4\rho^2}\mathfrak{B}^{*}(I_3+\mathfrak{A}+\Omega')\Omega$, and also
\begin{align*}
 \pPi_{1}Q_{2}&=\begin{bmatrix}
       \Omega&0&0\\
       -(I_3-\Omega)\mathfrak{A}\Omega-\Omega'\Omega&0&0\\
       0&0&0
      \end{bmatrix},\; 
      \pPi_{2}=\begin{bmatrix}
       I_3-\Omega&0&0\\
       (I_3-\Omega)\mathfrak{A}\Omega+\Omega'\Omega&I_3-\Omega&0\\
       0&0&0
      \end{bmatrix}.
\end{align*}
Regarding that $\im G_{0}=\im G_{1}=\im G_{2}$ and owing to \cite[Proposition 3.20]{CRR}  we determine the rank of $G_{3}=G_{2}+ B_{2} Q_{2}$ without knowing $B_{2}$ in detail:
\[\rank G_{3}=\rank G_{2}+\rank(I_7-G_{0}G_{0}^{+})B Q_{2}=6+1=7.
\]
The projector function\footnote{We call attention to an error in the projector representation in \cite[Section 6]{HM2020}. There the term $(I_3-\Omega)\mathfrak{A}\Omega$ is missing. Instead of the above correct $\pPi_{2}$ there the incorrect version
 \begin{align*}
 \begin{bmatrix}
       I_3-\Omega&0&0\\
       \Omega'\Omega&I_3-\Omega&0\\
       0&0&0
      \end{bmatrix} =:\tilde{\pPi}_{2}.
 \end{align*}
is given.}  $\pPi_{2}$ is not identical with the canonical projector function $\pPi_{can}$  since its image does not coincide with $S_{can}$. The determination of $\pPi_{can}$ according to \cite{CRR} is much more complex.
But owing to \cite[Theorem 2.8]{CRR} which lists invariances of the construction of admissible sequences, it holds that
\begin{align}\label{Ncan}
 N_{can}&=\ker \pPi_{can}=\ker \pPi_{2}\nonumber\\
 &=\{\begin{bmatrix}
      Z_{1}\\Z_{2}\\Z_{3}
     \end{bmatrix}\in \Real^3\times\Real^3\times\Real:Z_{1},Z_{2b}\in \im\mathfrak{B}, Z_{2}=Z_{2b}-(I_3-\Omega)\mathfrak{A}Z_{1}-\Omega'Z_{1}
\}\\
&=\im C_{N_{can}},\quad  C_{N_{can}}=\begin{bmatrix}
       \mathfrak B&0&0\\
       -((I_3-\Omega)\mathfrak A+\Omega')\mathfrak B&\mathfrak B&0\\
       0&0&1
      \end{bmatrix}.\nonumber
\end{align}

\subsection{Accurately stated initial conditions to \eqref{DAE1} by using \protect$\pPi_2\protect$}\label{subs:IC}

The initial condition to  \eqref{DAE1},
\begin{align}\label{IC}
 G_{a}x(a)=\gamma,
\end{align}
with a matrix $G_{a}\in \Real^{d\times m}$,
 is accurately stated, if 
\begin{align*}
 \im G_{a}=\Real^{d},\quad \ker G_{a}=\ker \pPi_{can}.
\end{align*}
We now intend to build a suitable matrix $G_{a}$. More precisely, we aim for a matrix function $G:\mathcal I\rightarrow\Real^{d\times m}$, such that $G(a)$ may serve as $G_{a}$ for each arbitrary $a\in \mathcal I$.

Introducing a matrix function $H:\mathcal I\rightarrow\Real^{2\times 3}$, with properties
\begin{align*}
 \rank H(t)=2,\quad \ker H(t)=\im\mathfrak{B}(t)=\im \Omega(t), \quad t\in \mathcal I,
\end{align*}
we form
\begin{align*}
 G=\begin{bmatrix}
    H&0&0\\0&H&0
   \end{bmatrix} \pPi_{2}=
   \begin{bmatrix}
    H&0&0\\H(\mathfrak{A}+\Omega')\Omega&H&0
   \end{bmatrix}.
\end{align*}
By construction, it holds that 
\begin{align*}
 \im G(t)=\Real^{d},\quad \ker G(t)=\ker \pPi_{2}(t)=\ker \pPi_{can}(t), \quad t\in\mathcal I.
\end{align*}
Below we choose\footnote{
Note that $H=V^*$ can be chosen, 
with any basis $V$ of $\ker \mathfrak{B}^*$.}
\begin{align*}
 H(t)=\begin{bmatrix}
       \sin t&-\cos t&0\\0&1&\cos t
      \end{bmatrix}.
\end{align*}
In particular, for $a=0$, this yields
\begin{align*}
 G(0)=\begin{bmatrix}
       0&-1&0&0&0&0&0\\
       0&1&1&0&0&0&0\\
       0&0&0&0&-1&0&0\\
       -1&0&0&0&1&1&0
      \end{bmatrix}.
\end{align*}
\subsection{Providing a matrix $G_a$ such that $\ker G_a=N_{can}(a)$ via a basis of $S_{*can}(a)$}\label{subs:adjoint}

We  write the homogeneous adjoint DAE to DAE \eqref{DAE1} in the form 
\begin{align}\label{adjDAE}
 \begin{bmatrix}
  X'_{1}\\
   X'_{2}\\
   0
 \end{bmatrix}+
 \begin{bmatrix}
       0&-\mathfrak{C}^{*}&\mathfrak{B}\\
       I&-\mathfrak{A}^{*}&0\\
       0&-\mathfrak{B}^*&0
      \end{bmatrix} \begin{bmatrix}
  X_{1}\\
   X_{2}\\
   X_{3}
 \end{bmatrix}= 0, \quad X_{1}=\begin{bmatrix}
  x_{1}\\
   x_{2}\\
   x_{3}
 \end{bmatrix},\; X_{2}=\begin{bmatrix}
  x_{4}\\
   x_{5}\\
   x_{6}
 \end{bmatrix},\; X_{3}=x_{7},
\end{align}
and start the reduction of the pair of  matrix functions featuring size $7\times7$, $\rank E_0=6$,
\begin{align*}
 E_0=\begin{bmatrix}
    I&0&0\\
    0&I&0\\
    0&0&0
   \end{bmatrix},\quad F_0=
   \begin{bmatrix}
       0&-\mathfrak{C}^{*}&\mathfrak{B}\\
       I&-\mathfrak{A}^{*}&0\\
       0&-\mathfrak{B}^*&0
      \end{bmatrix},
\end{align*}
by choosing 
\begin{align*}
 Y^*_0=\begin{bmatrix}
      I&0&0\\0&I&0
     \end{bmatrix},\quad
Z^*_0=\begin{bmatrix}
      0&0&1
     \end{bmatrix},
\end{align*}
and forming
\begin{align*}
 S_{0}=\ker Z^*F_0=\{\begin{bmatrix}
                        Z_1\\Z_2\\Z_3
                       \end{bmatrix}\in\Real^3\times\Real^3\times\Real:\mathfrak{B}^*Z_2=0
\},\quad \dim S_{0}=6.
\end{align*}
This leads to the following basis $C_{0}$ of the subspace $S_{0}$,
\begin{align*}
 C_{0}=\begin{bmatrix}
        I_3&0&0\\
        0&C_{\mathfrak{B}}&0\\
        \underbrace{0}_{3}&\underbrace{0}_{2}&\underbrace{1}_{1}
       \end{bmatrix},
\end{align*}
and also to the reduced pair (size $6\times 6$),
\begin{align*}
 &E_1=Y^*_0 EC_0=\begin{bmatrix}
    I_3&0&0\\
    0&C_{\mathfrak{B}}&0\\
   \end{bmatrix},\quad \rank E_1=5,\\
 &F_1=Y^*_0 FC_0+Y^*_0 EC'_0 =
   \begin{bmatrix}
       0&-\mathfrak{C}^{*}C_{\mathfrak{B}}&\mathfrak{B}\\
       I_3&-\mathfrak{A}^{*}C_{\mathfrak{B}}+C'_{\mathfrak{B}}&0
      \end{bmatrix},
\end{align*}
Remember that $C_{\mathfrak{B}}$ is a smooth basis of $\ker\mathfrak{B}^*$ \eqref{xx} and $\Omega$ is the orthoprojector onto $\im\mathfrak{B}$.
With 
\begin{align*}
 Y^*_1=\begin{bmatrix}
      I_3&0\\0&C^*_{\mathfrak{B}}
     \end{bmatrix},\quad
Z^*_1=\begin{bmatrix}
      0&\mathfrak{B}^*
     \end{bmatrix},
\end{align*}
and 
\begin{align*}
 &S_{1}=\ker Z^*_1F_1=\{\begin{bmatrix}
                        Z_1\\Z_2\\Z_3
                       \end{bmatrix}\in\Real^3\times\Real^2\times\Real:\mathfrak{B}^*Z_1+\mathfrak{B}^*(C'_{\mathfrak{B}}-\mathfrak{A}^*C_{\mathfrak{B}})Z_2=0
\},\\
&\dim S_{1}=5,
\end{align*}
we arrive at the smooth basis $C_1$ to the subspace $S_{1}$,
\begin{align*}
C_{1}=\begin{bmatrix}
        C_{\mathfrak{B}}&-\Omega(C'_{\mathfrak{B}}-\mathfrak{A}^*C_{\mathfrak{B}})&0\\
        0&I_2&0\\
        \underbrace{0}_{2}&\underbrace{0}_{2}&\underbrace{1}_{1}
       \end{bmatrix},
\end{align*}
as well as the next reduced pair, (size $5\times 5$),
\begin{align*}
 &E_2=Y^*_1 E_1C_1=\begin{bmatrix}
    C_{\mathfrak{B}}&-\mathfrak{F}&0\\
    0&C^*_{\mathfrak{B}}C_{\mathfrak{B}}&0\\
   \end{bmatrix},\quad 
   \mathfrak{F}=\Omega(C'_{\mathfrak{B}}-\mathfrak{A}^*C_{\mathfrak{B}}) \quad\rank E_2=4,\\
 &F_2=Y^*_1F_1C_1+Y^*_1E_1C'_1 =
   \begin{bmatrix}
       C'_{\mathfrak{B}}&-\mathfrak{C}^{*}C_{\mathfrak{B}}-\mathfrak{F}'&\mathfrak{B}\\
       C^*_{\mathfrak{B}}C_{\mathfrak{B}}&C^*_{\mathfrak{B}}(C_{\mathfrak{B}}'-\mathfrak{A}^{*}C_{\mathfrak{B}})&0
      \end{bmatrix},
\end{align*}
With 
\begin{align*}
Z^*_2=\begin{bmatrix}
      \mathfrak{B}^*&\;\mathfrak{B}^*\mathfrak{E}
     \end{bmatrix},\quad \mathfrak{E}=\Omega \mathfrak{F}(C^*_{\mathfrak{B}}C_{\mathfrak{B}})^{-1},
\end{align*}
we obtain
\begin{align*}
 &S_{2}=\ker Z^*_2F_2=\{\begin{bmatrix}
                        Z_1\\Z_2\\Z_3
                       \end{bmatrix}\in\Real^2\times\Real^2\times\Real:-\mathfrak{H}_{1}Z_1-\mathfrak{H}_{2}Z_2+Z_3=0
\},\\
&\mathfrak{H}_1=-\frac{1}{4\rho^2}\mathfrak{B}^*( \Omega C'_{\mathfrak{B}}+\mathfrak{F} ) ,\\ 
&\mathfrak{H}_2=-\frac{1}{4\rho^2}\mathfrak{B}^*(-\Omega(\mathfrak{C}^* C_{\mathfrak{B}}+\mathfrak{F}') +\mathfrak{E}C^*_{\mathfrak{B}}(C'_{\mathfrak{B}}-\mathfrak{A}^*C_{\mathfrak{B}}) )\\
&\dim S_{2}=4,
\end{align*}
as well as the basis
\begin{align*}
C_{2}=\begin{bmatrix}
        I_2&0\\0&I_2\\
        \underbrace{\mathfrak{H}_1}_{2}&\underbrace{\mathfrak{H}_2}_{2}
       \end{bmatrix},
\end{align*}
Finally, taking
\begin{align*}
Y^*_2=\begin{bmatrix}
      C_{\mathfrak{B}}^*&0\\
      0&I_2
     \end{bmatrix},
\end{align*}
the next matrix function $E_3$ (size $4\times 4$) remains nonsingular, namely
\begin{align*}
 E_3=Y^*_2 E_2C_2=\begin{bmatrix}
   C_{\mathfrak{B}}^* C_{\mathfrak{B}}&-C_{\mathfrak{B}}^*\mathfrak{F}\\
    0&C^*_{\mathfrak{B}}C_{\mathfrak{B}}\\
   \end{bmatrix},
\end{align*}
and, hence the matrix function $C=C_0C_1C_2$ showing size $7\times 4$,
\begin{align*}
 C&=\begin{bmatrix}
        I_3&0&0\\
        0&C_{\mathfrak{B}}&0\\
        \underbrace{0}_{3}&\underbrace{0}_{2}&\underbrace{1}_{1}
       \end{bmatrix}
    \begin{bmatrix}
        C_{\mathfrak{B}}&-\Omega(C'_{\mathfrak{B}}-\mathfrak{A}^*C_{\mathfrak{B}})&0\\
        0&I_2&0\\
        \underbrace{0}_{2}&\underbrace{0}_{2}&\underbrace{1}_{1}
       \end{bmatrix}
     \begin{bmatrix}
        I_2&0\\0&I_2\\
        \underbrace{\mathfrak{H}_1}_{2}&\underbrace{\mathfrak{H}_2}_{2}
       \end{bmatrix}\\
       &=\begin{bmatrix}
          C_{\mathfrak{B}}&-\Omega(C'_{\mathfrak{B}}-\mathfrak{A}^*C_{\mathfrak{B}})\\
          0&C_{\mathfrak{B}}\\
           \mathfrak{H}_1&\mathfrak{H}_2
         \end{bmatrix}
\end{align*}
serves as basis of $S_{*\,can}$\footnote{Note that \eqref{adjDAE} is the adjoint to \eqref{DAE1}.}.
We are mainly interested in the matrix function
\begin{align*}
 C^*E_0&= \begin{bmatrix}
          C_{\mathfrak{B}}^*&0&\mathfrak{H}^*_1\\
          (-\Omega(C'_{\mathfrak{B}}-\mathfrak{A}^*C_{\mathfrak{B}}))^*&
          C_{\mathfrak{B}}^*&
          \mathfrak{H}^*_2
         \end{bmatrix}
         \begin{bmatrix}
          I&0&0\\0&I&0\\
          0&0&0
         \end{bmatrix}\\
         &=
         \begin{bmatrix}
          C_{\mathfrak{B}}^*&0&0\\
          (-\Omega(C'_{\mathfrak{B}}-\mathfrak{A}^*C_{\mathfrak{B}}))^*&
          C_{\mathfrak{B}}^*&0
         \end{bmatrix}\\
       &=\begin{bmatrix}
          C_{\mathfrak{B}}^*&0&0\\
          ((\Omega'+\mathfrak{A}^*)C_{\mathfrak{B}}))^*&
          C_{\mathfrak{B}}^*&0
         \end{bmatrix} 
         =\begin{bmatrix}
          C_{\mathfrak{B}}^*&0&0\\
          C_{\mathfrak{B}}^*(\Omega'+\mathfrak{A})&
          C_{\mathfrak{B}}^*&
          0
         \end{bmatrix}  
\end{align*}
and its three-dimensional nullspace
\begin{align*}
 \ker C^*E=\{\begin{bmatrix}
                        Z_1\\Z_2\\Z_3
                       \end{bmatrix}\in\Real^3\times\Real^3\times\Real:
                       Z_2=Z_{2b}-(I-\Omega)(\Omega'+\mathfrak{A})Z_1, \;
                       Z_1, Z_{2b}\in \im\mathfrak{B}
\}
\end{align*}
which coincides with the subspace $N_{can}$ in \eqref{Ncan}.
\subsection{Canonical projector function to \eqref{DAE1}}
Having the bases functions  $C_{S_{can}}$, $C_{N_{can}}$ of $S_{can}$, $N_{can}$, respectively, the $7\times 7$ matrix function $\mathcal M:=[C_{S_{can}}\,C_{N_{can}}]$  remains nonsingular everywhere on $\mathcal I$, and 
\begin{align*}
 \pPi_{can}=\mathcal M
 \begin{bmatrix}
             I_4&\\&0_3
            \end{bmatrix} \mathcal M^{-1}.
\end{align*}

\appendix
\section{Appendix}\label{s.Appendix}
\subsection{Standard form DAEs, DAEs with properly involved derivative, and proper factorizations}\label{subs.forms}
We say that \emph{$N$ is a $C^{1}$- subspace in $\Real^{n}$}, if  $N(t)\subseteq \Real^{n}$ is a time-varying subspace, $t\in \mathcal I$, and the projector-valued function $Q:\mathcal I\rightarrow \Real^{n\times n}$, with $Q(t)=Q(t)^{2}=Q(t)^{T}$, $\im Q(t)=N(t)$, $t\in\mathcal I$, is continuously differentiable. We underline  that any $C^{1}$- subspace in $\Real^{n}$ has constant dimension and also a continuously differentiable basis.

Each continuous  matrix function $E:\mathcal I\rightarrow\Real^{m\times m}$ having  constant rank $r$ and a nullspace which is a $C^{1}$-subspace in $\Real^{m}$ can be factorized into $E=AD$, with  $A:\mathcal I\rightarrow\Real^{m\times n},\; D:\mathcal I\rightarrow \Real^{n\times m}$, so that $A$ is continuous, $D$ is continuously differentiable, $\ker A$ and $\im D$ are a $C^{1}$-subspaces in $\Real^{n}$, and 
\begin{align}
 \ker A(t)\oplus \im D(t)&=\Real^{n},\quad t\in\mathcal I.\label{transver}
\end{align}
A possible choice is $n=m$,\; $A=E,\,D=E^{+}E$. If $E$ is itself continuously differentiable, then also the factor $A$ can be chosen to be continuously differentiable, for instance $A=EE^{+}$, $D=E$.
Any factorization satisfying  condition \eqref{transver} is called \emph{proper factorization}. Note that then the function $R:\mathcal I\rightarrow \Real^{n\times n}$, projecting pointwise onto $\im D$ along $\ker A$ is also continuously differentiable and one has
\begin{align}\label{R}
 \im E=\im A,\;\ker E=\ker D,\quad A=AR,\;D=RD.
\end{align}
$R$ is then called \emph{border-projector function}.
\medskip

Using any proper factorization of the leading coefficient $E$, the standard form DAE 
$Ex'+Fx=q$
can be rewritten with $B=F-AD'$ as \emph{DAE with properly stated leading term or DAE with properly involved derivative},
$
 A(Dx)'+Bx=q.
$. 

Of course, on the other hand, starting from a DAE with properly involved derivative, $A(Dx)'+Bx=q$, one immediately gains the DAE in standard form \linebreak
$ADx'+(B+AD')x=q.$

A properly involved derivative is essential when rigorous solvability statements are required and only the component $Dx$ can be expected to be continuously differentiable. 
In the case of DAEs with proper involved derivative, the DAE and its adjoint show a formal symmetry that proves beneficial for the analysis.
However, if, as here, we are concerned with the description of the time-varying subspaces in $\Real^m$ that capture the solution values, then the form of the inclusion of the derivative does not matter, in particular,
\begin{align*}
 N&=\ker E=\ker AD=\ker D,\\
 S&=\{z\in\Real^m:Fz\in \im E\}=\{z\in\Real^m:Bz\in \im AD\},
\end{align*}
and  both DAE forms share their canonical projector $\pPi_{can}$ and the related canonical subspaces $S_{can}$ and $N_{can}$.

\subsection{Modification of the reduction procedure for DAEs  \eqref{pDAE}}\label{subs.proper}
We sketch here a modification of the reduction procedure from \cite{RaRh} for DAEs
with properly involved derivative \eqref{pDAE}.

We start by $A_0=A,\, D_0=D,\,B_0=B,\,m_0=m,\,r_{0}=r$ and consider the homogeneous DAE
\begin{align*}
 A_0(D_0x)'+B_0=0.
\end{align*}
By means of a basis $Z_0:\mathcal I\rightarrow \Real^{(m_0-r_0)\times m_0}$ of $(\im A_0)^{\perp}=\ker A_0^{*}$ and a basis $Y_0:\mathcal I\rightarrow \Real^{r_0\times m_0}$ of $\im A_0$ we divide the DAE into the two parts
\begin{align*}
 Y_0^*A_0(D_0x)'+Y_0^*B_0x=0,\quad Z_0^*B_0x=0.
\end{align*}

From  $\im[A_0D_0\;\;B_0]=\Real^m$ we derive that $\rank Z_0^*B= m_0-r_{0}$, and hence. the subspace   
$S_{0}=\ker Z_0^*B=\ker Z_0^*F$ has dimension $r_{0}$. Each solution of the homogeneous DAE must stay in the subspace $S_{0}$. Choosing a basis 
 $C_0:\mathcal I\rightarrow \Real^{r_0\times m_0}$ of $S_{0}$, each solution of the DAE can be represented as $x=C_0 x_{[1]}$, with a function $x_{[1]}:\mathcal I\rightarrow \Real^{r_0}$ satisfying the reduced to size $m_1=r_{0}$ DAE given below.
 In contrast to \cite{RaRh} where the basis $C_0$ is required to be continuously differentiable, we suppose now a continuous basis $C_0$ which has a continuosly differentiable part $D_0C_0$. Using a pointwise generalized inverse $(D_0C_0)^-$ of $D_0C_0$ we may write
 \begin{align*}
  D_0x&=D_0C_0x_{[1]}=D_0C_0\underbrace{(D_0C_0)^-D_0C_0}_{=:D_1}x_{[1]},\\
  (D_0x)'&=D_0C_0(D_1x_{[1]})'+(D_0C_0)'D_1x_{[1]}.
 \end{align*}
This leads to a DAE living in $\Real^{m_1}$, $m_1:=r_{0}$, with properly involved derivative, 
 \begin{align*}
 \underbrace{Y_0^*A_0D_0C_0}_{=:A_1}(D_1 x_{[1]})'+\underbrace{Y_0^*(B_0C_0+A_0(D_0C_0)'D_1)}_{=:B_1}x_{[1]}=0.
\end{align*}
\subsection{Basic steps by Cistyakov and Jansen}\label{subs.CistJans}
Let the pair $\{E,F\}$, $E,F:\mathcal I\rightarrow\Real^{m\times m}$, be pre-regular with constants $r$ and $\theta$ according to Definition \ref{d.prereg}.
We take over some notations from Section \ref{s.comparison}. 

Let $T,T^c,Z$, and $Y $ represent bases of $\ker E, (\ker E)^{\perp}, (\im E)^{\perp}$, and $\im E $, respectively.

The matrix function $Z^*FT$ has size $(m-r)\times(m-r)$, rank $m-r-\theta=a$, and $\ker Z^*FT =T^+(N\cap S)$ has dimension $\theta$.
By scaling with $ [Y\, Z]^*$ one splits  the DAE
\begin{align*}
 Ex'+Fx=q
\end{align*}
 into the partitioned shape
\begin{align}
 Y^*Ex'+Y^*Fx&=Y^*q,\label{A.1}\\
 Z^*Fx&=Z^*q.\label{A.2}
\end{align}
Owing to the pre-regularity, the $(m-r)\times m$ matrix function  $ Z^*F$ features full row-rank $m-r$. We keep in mind that $S=\ker Z^*F$ has dimension $r$.
\subsubsection{Elimination by Cistyakov}
Taking a nonsingular matrix function  $K$ of size $m\times m$ such that $Z^*FK=: [\tilde F_{21}\, \tilde F_{22}]$, with $\tilde F_{22}$ being nonsingular, the transformation  $x=K\tilde x$ turns  \eqref{A.2} into 
\begin{align*}
Z^*Fx&=Z^*FK\tilde x=:\tilde F_{21}\tilde x_{1}+\tilde F_{22}\tilde x_{2}=Z^*q, \\
&\text{yielding}\quad \tilde x_2=-\tilde F_{22}^{-1}\tilde F_{21}\tilde x_1+Z^*q.
\end{align*}
Next one eliminates the variable $\tilde x_2$ in the   transformed version of \eqref{A.1},
\begin{align*}
 \underbrace{Y^*EK}_{= :[\tilde E_{11}\, \tilde E_{12}]} \tilde x'+(Y^*FK + Y^*EK')\tilde x=Y^*q,
\end{align*}
which yields a DAE for $\tilde x_1=:x_{[1]}$,
\begin{align*}
 \underbrace{(\tilde E_{11}-\tilde E_{12}\tilde F_{22}^{-1}\tilde F_{21})}_{=: E_{[1]}} \tilde x_1'+ F_{[1]} \tilde x_1=Y^*q.
\end{align*}
A further look at the matter shows that we are dealing with a special basis of the subspace $S$, namely
\begin{align*}
 S=\im C,\quad C= K \begin{bmatrix}
                  I_r\\ -\tilde F_{22}^{-1}\tilde F_{21}
                 \end{bmatrix},
\end{align*}
whereby this early predecessor of the method described in \cite{RaRh} can now be classified as its special version. Note that the procedure in \cite{RaRh} allows for the choice of an arbitrary basis for $S$.

It should be further mentioned, that in \cite{Cis1982} the elimination method is applied not only for pre-regular pairs but for general rectangular matrix functions $E,F$.
\subsubsection{Dissectionn by Jansen}
The approach in \cite{Jansen2014} needs several more splittings. As before let
$T,T^c, Z$, and $Y$ represent bases of $\ker E, (\ker E)^{\perp}, (\im E)^{\perp}$, and $\im E$, respectively. Additionally, let $V,W$ be bases of  $(\im Z^*FT)^{\perp}$, and $\im Z^*FT$. By construction, $V$ has size $(m-r)\times a$ and $W$ has size $(m-r)\times \theta$.
One  starts with the  transformation
\begin{align*}
 x= \begin{bmatrix}
                       T^c& T
                      \end{bmatrix}\tilde x, \quad 
                      \tilde x=\begin{bmatrix}
                       \tilde x_1\\\tilde x_2
                      \end{bmatrix},\quad x=T^c\tilde x_1+ T\tilde x_2.
\end{align*}
The background is the associated possibility to suppress the derivative of the nullspace-part $T\tilde x_n$ similarly as in the context of properly formulated DAEs and to set $Ex'= ET^c\tilde x_1'+E{T^c}'\tilde x_1 +ET'\tilde x_2$, which, however, does not play a role here where altogether continuously differentiable solutions are assumed. Furthermore, an additional partition of the derivative-free equation \eqref{A.2} by means of the scaling with $[V\,W]^*$ is applied, which results in the system
\begin{align}
Y^*ET^c\tilde x'_1 +Y^*(FT^c+E{T^c}')\tilde x_1+ Y^*(FT+E{T}')\tilde x_2&=Y^*q,\label{A3}\\
V^*Z^*FT^c\tilde x_1+V^*Z^*FT\tilde x_2&=V^*Z^*q,\label{A4}\\
W^*Z^*FT^c\tilde x_1 \hspace*{18mm} &=W^*Z^*q,\label{A5}.
\end{align}
The matrix function $W^*Z^*FT^c$ has full row-rank $\theta$ and $V^*Z^*FT$ has full row-rank $a$. Now comes another split. Choosing bases $G, H$ of $\ker W^*Z^*FT^c\subset\Real^{\theta}$ and $\ker V^*Z^*FT\subset\Real^{a}$, as well as bases of respective complementary subspaces, we transform 
\begin{align*}
 \tilde x_1= \begin{bmatrix}
                       G^c& G
                      \end{bmatrix}\bar x_1\quad 
                      \bar x_1=\begin{bmatrix}
                       \bar x_{1,1}\\\bar x_{1,2}
                      \end{bmatrix},\quad \tilde x_1=G^c\bar x_{1,1}+ G\bar x_{1,2},\\
 \tilde x_2= \begin{bmatrix}
                       H^c& H
                      \end{bmatrix}\bar x_2\quad 
                      \bar x_2=\begin{bmatrix}
                       \bar x_{2,1}\\\bar x_{2,2}
                      \end{bmatrix},\quad \tilde x_2=H^c\bar x_{2,1}+ H\bar x_{2,2}.                     
\end{align*}
Thus equations \eqref{A4} and \eqref{A5} are split into
\begin{align*}
 V^*Z^*FT^c(G^c\bar x_{1,1}+ G\bar x_{1,2})+V^*Z^*FT H^c\bar x_{2,1}&=V^*Z^*q,\\
W^*Z^*FT^c G^c\bar x_{1,1} \hspace*{38mm} &=W^*Z^*q.
\end{align*}
The matrix functions $V^*Z^*FT H^c$ and $W^*Z^*FT^c G^c$ are nonsingular each, which allows the resolution to $\bar x_{1,1}$ and $\bar x_{2,1}$. In particular, for $q=0$ it results that 
$\bar x_{1,1}= 0$ and $\bar x_{2,1}= \mathfrak E \bar x_{1,2}$, with 
\[\mathfrak E:=-(V^*Z^*FT H^c)^{-1}V^*Z^*FT^cG.
\]
Overall, therefore, the latter procedure presents again a transformation, namely
\begin{align*}
 x= K\bar x,\quad K=\begin{bmatrix}
    T^cG^c&T^cG&TH^c&TH
   \end{bmatrix},\quad \bar x=\begin{bmatrix}
   \bar x_{1,1}\\\bar x_{1,2}\\\bar x_{2,1}\\\bar x_{2,2}
   \end{bmatrix}
   \in \Real^{\theta}\times\Real^{r-\theta}\times\Real^{a}\times\Real^{\theta},
\end{align*}
and we realize that we have found again a basis of the subspace $S$, namely
\begin{align*}
 S=\im C,\quad C=K
 \begin{bmatrix}
    0&0\\I_{r-\theta}&0\\\mathfrak E&0\\0&I_{\theta}
                  \end{bmatrix}=
\begin{bmatrix}
 T^cG+TH^c\mathfrak E &\; TH
\end{bmatrix},
\end{align*}
 which makes the dissection approach a special case of \cite{RaRh}. 
 The characteristic values together with the index are formally adapted to the values of the tractability index. It starts with $r^{D}_0=r$, and is continued in ascending order with $r^{D}_{i+1}=r^{D}_{i}+ a_{i}=r^{D}_{i}+\rank Z^+_{i}F_iT_i$ etc. until $r^{D}_{\mu-1}<r^{D}_{\mu}=m$, \cite[Definition 4.13]{Jansen2014}.

It should be noted, however, that the dissection concept is developed with considerable  effort for nonlinear DAEs in \cite{Jansen2014}. 
\subsection{Admissible matrix function sequences and related subspaces}\label{AdmissibleMatrix}

In this part we apply several routine notations and tools used in the projector based analysis of DAEs. We refer to the appendix for a short roundup and to \cite{CRR,Mae2014} for more details.

Given are at least continuous matrix functions $E,F:\mathcal I\rightarrow\Real^{m\times m}$, $E$ has a $C^{1}$-nullspace and constant rank $r$. We use a proper factorization $E=AD$ where $A:\mathcal I\rightarrow\Real^{m\times k}$, $D:\mathcal I\rightarrow\Real^{k\times m}$, and $B=-(F+AD')$.
$R:\mathcal I:\rightarrow\Real^{k\times k}$ denotes the continuously differentiable projector-valued function such that $\im D=\im R$ and $\ker A=\ker R$.

Let $Q_{0}:\mathcal I\rightarrow \Real^{m\times m}$ denote any continuously differentiable projector-valued function such that $\im Q_{0}=\ker D=\ker E$, for instance, $Q_{0}=I-D^{+}D$ with the pointwise Moore-Penrose inverse $D^+$. Set $P_{0}=I-Q_{0}$ and 
let $D^{-}$ denote the pointwise generalized inverse of $D$ determined by
\[
 D^{-}DD^{-}=D^{-},\quad DD^{-}D=D,\quad DD^{-}=R,\quad D^{-}D=P_{0}.
\]
Set $G_{0}=AD,\; B_{0}=B,\; \pPi_{0}=P_{0}$. For a given level $\kappa\in \Natu$, the sequence $G_{0},\ldots, G_{\kappa}$ is called an \emph{admissible matrix function sequence} associated with the pair $\{E,F\}$ and triple $\{A,D,B\}$, respectively, e.g.,\cite[Definition 2.6]{CRR} if it is built by the rule
\begin{align*}
 G_{i}=&G_{i-1}+B_{i-1}Q_{i-1},\\
 &B_{i}=B_{i-1}P_{i-1}-G_{i}D^{-}(D\pPi_{i}D^{-})'D\pPi_{i-1},\\
 &N_{i}=\ker G_{i}, \quad \widehat{N_{i}}:=(N_{0}+\cdots+N_{i-1})\cap N_{i}, \quad N_{0}+\cdots+N_{i-1}=:\widehat{N_{i}}\oplus X_{i},\\
 &\text{choose } Q_{i} \text{ such that } Q_{i}=Q_{i}^{2}, \; \im Q_{i}=N_{i},\; X_{i}\subseteq \ker Q_{i},\\
 &P_{i}=I-Q_{i},\; \pPi_{i}=\pPi_{i-1}P_{i},\\
 i=1&,\ldots,\kappa,
\end{align*}
and, additionally, 
\begin{description}
 \item[\textrm{(a)}] $G_{i}$ has constant rank $r^T_{i}$,\; $i=0,\ldots,\kappa$,
 \item[\textrm{(b)}] $\widehat{N_{i}}=(\im [\pPi_{i-1}^{*}\,G_{i}^{*}])^{\perp}$ has constant dimension $u^T_{i}$,\; $i=1,\ldots,\kappa$,
 \item[\textrm{(c)}] $\pPi_{i}$ is continuous and $D\pPi_{i}D^{-}$ is continuously differentiable,\; $i=0,\ldots,\kappa$.
\end{description}
The admissible matrix functions $G_{i}$ are continuous. The construction is supported by two constant-rank conditions at each level. It results that
\begin{align*}
 0<r^T_{0}\leq r^T_{1}\leq\cdots\leq r^T_{i}\leq \ldots,\\
 0\leq u^T_{1}\leq u^T_{2}\leq\cdots\leq u^T_{i}\leq \ldots
\end{align*}
By construction, the inclusions
\begin{align}\label{Gsequence}
 \im G_{0}\subseteq \im G_{1}\subseteq\ldots\subseteq\im G_{r+1}=\im G_{r+2}
\end{align}
are valid pointwise.
There are several special projector functions incorporated in an admissible matrix function sequence, among them admissible projectors $Q_{i}$ onto $\ker G_{i}$ and $\pPi_{i}=\pPi_{i-1}(I-Q_{i})$, $\pPi_{0}=(I-Q_{0})$, yielding the further inclusions
\begin{align}\label{Psequence}
 \ker \pPi_{0}\subseteq \ker \pPi_{1}\subseteq\ldots\subseteq\ker \pPi_{r}=\ker \pPi_{r+1}.
\end{align}
Each of the time-varying subspaces in \eqref{Gsequence} and \eqref{Psequence} has constant dimension, which is ensured by the respective rank conditions.

The subspaces involved in \eqref{Gsequence} and \eqref{Psequence} are proved to be invariant with respect to special possible choices within the construction procedure and also with respect to the factorization of $E=AD$.


\begin{definition}\label{d.tractability}
The DAE given by the  coefficient function pair  $\{E,F\}$ or a related tripel $\{A,D,B\}$ is called \emph{regular} if there are an index $\mu$ and an  admissible matrix function sequence $G_0,G_1,\ldots,
G_{\mu}$ such that $r^T_{\mu}=m$.
The \emph{tractability index} of the DAE is defined to be the smallest index  $\mu$ with $r^T_{\mu}=m$, and 
the integers 
\[0<r^T_{0}\leq r^T_{1}\leq\cdots\leq r^T_{\mu-1} <r^T_{\mu}=m
\]
are called \emph{characteristic values} of the regular DAE.
\end{definition}
For a regular DAE  $N_{\mu}=\{0\}$ is trivially valid, therefore  $\widehat{N_{\mu}}=\{0\}$ must also be valid
and in turn $u^T_{i}=0$, for all $i$. Moreover, owing to \cite[Proposition 2.5]{CRR} regularity requires 
\begin{align*}
 \im [E \; F]=\Real^{m}, \quad \im [AD \; B]=\Real^{m}.
\end{align*}

\bibliographystyle{plain}
\bibliography{hmInitialCond}

\begin{minipage}{0.4\textwidth}\authorlist
\end{minipage}

\end{document}